\documentclass[onefignum,onetabnum]{siamart190516}
\usepackage{amsfonts}
\usepackage{graphicx}
\usepackage{epstopdf}
\usepackage{algorithmic}
\ifpdf
  \DeclareGraphicsExtensions{.eps,.pdf,.png,.jpg}
\else
  \DeclareGraphicsExtensions{.eps}
\fi

\usepackage{siunitx}
\usepackage{listings}
\usepackage{graphicx,url}
\usepackage[utf8]{inputenc}
\usepackage{amsmath,amssymb,blkarray,epsfig,graphics,psfrag,latexsym,amsmath,listings,graphicx,url,tabularx,stmaryrd,booktabs}
\usepackage[utf8]{inputenc}
\usepackage{multicol}
\usepackage{bbm}
\usepackage{color}
\graphicspath{ {Figures/} }

\newcommand{\diam}{\mathrm{diam}}
\usepackage{color,soul}

\usepackage{tikz}
\usetikzlibrary{shapes}
\usetikzlibrary{matrix}
\usetikzlibrary{positioning}
\usetikzlibrary{fit}

\tikzstyle{res}=[circle,thick,minimum size=4mm,draw=black,fill=red,inner sep=1pt]
\tikzstyle{non-res}=[circle,thick,minimum size=4mm,draw=black,inner sep=1pt]
\tikzstyle{light-res}=[circle,thick,minimum size=4mm,draw=black,fill=red!40,inner sep=1pt]
\tikzstyle{blue}=[circle,thick,minimum size=4mm,draw=black,fill=blue!20,inner sep=1pt]

\newsiamremark{remark}{Remark}
\newsiamremark{hypothesis}{Hypothesis}
\crefname{hypothesis}{Hypothesis}{Hypotheses}
\newsiamthm{claim}{Claim}

\headers{A Survey of Metric Dimension and its Applications}{Richard C. Tillquist, Rafael M. Frongillo, Manuel E. Lladser}

\title{Getting the Lay of the Land in Discrete Space: A Survey of Metric Dimension and its Applications\thanks{Submitted to the editors 4/1/2021.
\funding{This article was partially funded by NSF IIS grant 1836914.}}}

\author{Richard C. Tillquist\thanks{Department of Computer Science, University of Colorado, Boulder, The United States (\email{richard.tillquist@colorado.edu}, \email{raf@colorado.edu})}
\and
Rafael M. Frongillo\footnotemark[2]
\and
Manuel E. Lladser\thanks{Department of Applied Mathematics, University of Colorado, Boulder, The United States (\email{manuel.lladser@colorado.edu})}
}

\usepackage{amsopn}

\ifpdf
\hypersetup{
  pdftitle={Getting the Lay of the Land in Discrete Space: A Survey of Metric Dimension and its Applications},
  pdfauthor={Richard C. Tillquist, Rafael M. Frongillo, and Manuel E. Lladser}
}
\fi

\begin{document}

\maketitle

\begin{abstract}
The metric dimension of a graph is the smallest number of nodes required to identify all other nodes based on shortest path distances uniquely. Applications of metric dimension include discovering the source of a spread in a network, canonically labeling graphs, and embedding symbolic data in low-dimensional Euclidean spaces. This survey gives a self-contained introduction to metric dimension and an overview of the quintessential results and applications. We discuss methods for approximating the metric dimension of general graphs, and specific bounds and asymptotic behavior for deterministic and random families of graphs. We conclude with related concepts and directions for future work.
\end{abstract}



\tableofcontents

\section{Introduction}

In the Euclidean plane, any set of three non-collinear points is enough to uniquely distinguish all points in the space based on distances. This process, called trilateration in $\mathbb{R}^2$, is the basic technique through which Global Positioning Systems (GPS) are able to pinpoint a location on the surface of the Earth. More generally, if $\|\cdot\|$ denotes the Euclidean distance and $R=\{\mathbf{r}_1,\ldots,\mathbf{r}_{n+1}\}\subset\mathbb{R}^n$ is a set of $n+1$ affinely independent points, the vectors $(\lVert \mathbf{x}-\mathbf{r}_1 \rVert, \dots, \lVert \mathbf{x}-\mathbf{r}_{n+1} \rVert)$ and $(\lVert \mathbf{y}-\mathbf{r}_1 \rVert, \dots, \lVert \mathbf{y}-\mathbf{r}_{n+1} \rVert)$ for $\mathbf{x},\mathbf{y} \in \mathbb{R}^n$ are different when $\mathbf{x} \neq \mathbf{y}$.

The situation becomes more complex, however, if the space of interest is discrete instead of continuous. One class of discrete spaces of particular interest are those which can be represented as graphs coupled with shortest path distance. On a graph $G=(V,E)$ the notion of \emph{metric dimension} is analogous to the number of satellites required for GPS to work effectively.
The goal is to pick a small set of vertices $R\subseteq V$ capable of identifying every vertex based solely on shortest path distances to $R$.
Solving this problem exactly is computationally complex but provides information useful in a variety of settings. A small set of ``satellites'' or ``landmarks'' in a discrete space can be valuable in assisting robots navigating over a physical space or in tracking the progress of a disease as it spreads between cities.
It could also be used in more abstract settings like identifying a source of misinformation in a social network, comparing network structure, categorizing chemical structures, or representing symbolic data numerically.

In this work we collate and interpret a number of theoretical results and approximation techniques associated with metric dimension, paying particular attention to specific types of graphs and applications.
We survey recent work and describe promising directions for future work.

\section{Formal Definition}
\label{sec:formal_definition}

Let $G=(V,E)$ be a graph, potentially with weighted edges, multi-edges, and self loops, and let $d(u,v)$ denote the shortest path distance in $G$ from $u \in V$ to $v \in V$.

\begin{definition}{(Resolving Set.)}
$R \subseteq V$ is resolving if, for all distinct $u, v \in V$, there exists $r \in R$ such that $d(r, u) \neq d(r, v)$. Such an $r$ is said to resolve or distinguish $u$ and $v$.
\end{definition}

\noindent In the context of $\mathbb{R}^n$, any set of $(n+1)$ or more affinely dependent points is analogous to a resolving set in a graph. 

By definition, $R=\{v_1,\ldots,v_k\}$ is resolving if and only if the transformation
\begin{displaymath}
d(u|R):=\big(d(r_1,u), \dots, d(r_k,u)\big)
\end{displaymath}
from $V$ to $\mathbb R^{|R|}$ is injective, i.e., every vertex $u \in V$ is uniquely represented by the vector of distances from all vertices in $R$ (listed in an arbitrary but specified order) to $u$. In many settings, minimizing the dimension $|R|$ of these vectors is a central goal.

\begin{definition}{(Metric Dimension.)}
The metric dimension $\beta(G)$ of $G$ is the smallest size of resolving sets on $G$. If $R$ is a resolving set on $G$ and $|R| = \beta(G)$, $R$ is called a minimal resolving, basis, or reference set of $G$.
\end{definition}

In the context of graphs, the concept of metric dimension was introduced separately by Slater in 1975~\cite{slater1975leaves} and by Harary and Melter in 1976~\cite{harary1976metric}, though the dimension of graphs was discussed earlier by Erd\"os et al. in 1965~\cite{erdos1965dimension}. Both the Slater and Harary and Melter papers focus on the metric dimension of trees and describe equivalent exact formula for graphs of this kind. 
Harary and Melter briefly discuss the metric dimension of several other types of graphs including cycles, complete graphs, and complete bipartite graphs though the metric dimension of wheel graphs is incorrectly stated as two. They also give an algorithm to reconstruct a tree given distances from every node to the elements of a resolving set. This is not possible for general graphs as not all edges are guaranteed to be represented in a shortest path with an element of a resolving set as an endpoint (\Cref{sec:strong-metric-dim}).

\subsection{Simple Examples}
\label{sec:simple_examples}

Before continuing, we examine several types of graphs for which minimal resolving sets are readily described and easily visualized. Through this examination, we hope to strengthen the reader's intuitive grasp of metric dimension and to solidify concepts that are important in future sections.

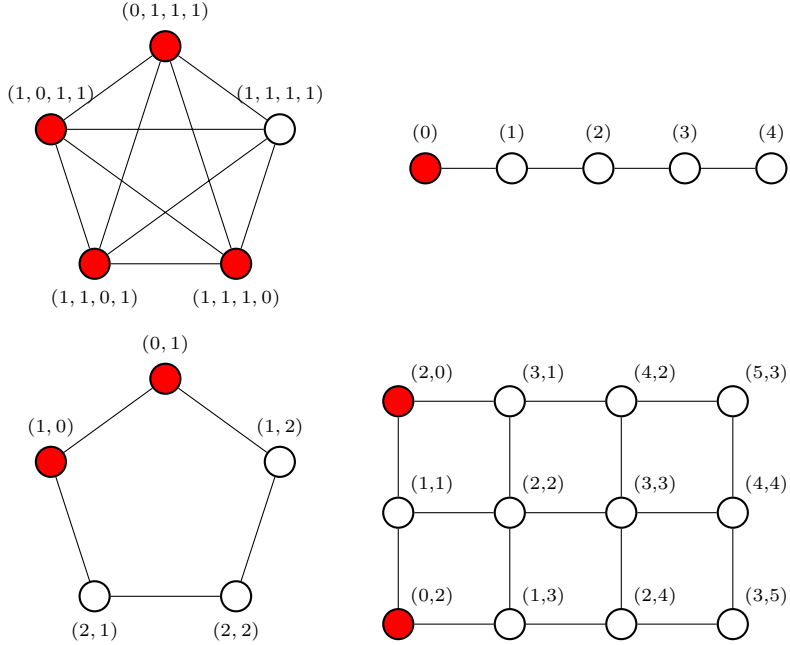
\begin{figure}[h]
\begin{tabular}{cc}
\begin{tikzpicture}[every node/.style={font=\scriptsize},scale=0.8]
    \def \n {5}
    \def \radius {2cm}
    \def \margin {8} 
    \node[res,label={$(0,1,1,1)$}] (A) at ({90+360/\n * (1 - 1)}:\radius) {};
    \node[res,label={$(1,0,1,1)$}] (B) at ({90+360/\n * (2 - 1)}:\radius) {};
    \node[res,label=below:{$(1,1,0,1)$}] (C) at ({90+360/\n * (3 - 1)}:\radius) {};
    \node[res,label=below:{$(1,1,1,0)$}] (D) at ({90+360/\n * (4 - 1)}:\radius) {};
    \node[non-res,label={$(1,1,1,1)$}] (E) at ({90+360/\n * (5 - 1)}:\radius) {};
    \draw (A) -- (B) -- (C) -- (D) -- (E) -- (A);
    \draw (A) -- (C) -- (E) -- (B) -- (D) -- (A);
\end{tikzpicture}
&
\raisebox{50pt}{~~~
\begin{tikzpicture}[every node/.style={font=\scriptsize},scale=1.15]
    \node[res,label=$(0)$] (A) at (0,0) {};
    \node[non-res,label=$(1)$] (B) at (1,0) {};
    \node[non-res,label=$(2)$] (C) at (2,0) {};
    \node[non-res,label=$(3)$] (D) at (3,0) {};
    \node[non-res,label=$(4)$] (E) at (4,0) {};
    \draw (A) -- (B) -- (C) -- (D) -- (E);
\end{tikzpicture}}
\\
\begin{tikzpicture}[every node/.style={font=\scriptsize},scale=0.8]
    \def \n {5}
    \def \radius {2cm}
    \def \margin {8} 
    \node[res,label={$(0,1)$}] (A) at ({90+360/\n * (1 - 1)}:\radius) {};
    \node[res,label={$(1,0)$}] (B) at ({90+360/\n * (2 - 1)}:\radius) {};
    \node[non-res,label=below:{$(2,1)$}] (C) at ({90+360/\n * (3 - 1)}:\radius) {};
    \node[non-res,label=below:{$(2,2)$}] (D) at ({90+360/\n * (4 - 1)}:\radius) {};
    \node[non-res,label={$(1,2)$}] (E) at ({90+360/\n * (5 - 1)}:\radius) {};
    \draw (A) -- (B) -- (C) -- (D) -- (E) -- (A);
\end{tikzpicture}
&
\begin{tikzpicture}[every node/.style={font=\scriptsize},scale=0.8]
  \matrix (m) [matrix of nodes,row sep=3em,column sep=3em,minimum width=2em,nodes={non-res}]
  {
     |[res]| & ~ & ~ & ~ \\
     ~ & ~ & ~ & ~ \\
     |[res]| & ~ & ~ & ~ \\
    };
  \draw (m-1-1) -- (m-1-2) -- (m-1-3) -- (m-1-4);
  \draw (m-2-1) -- (m-2-2) -- (m-2-3) -- (m-2-4);
  \draw (m-3-1) -- (m-3-2) -- (m-3-3) -- (m-3-4);
  
  \draw (m-1-1) -- (m-2-1) -- (m-3-1);
  \draw (m-1-2) -- (m-2-2) -- (m-3-2);
  \draw (m-1-3) -- (m-2-3) -- (m-3-3);
  \draw (m-1-4) -- (m-2-4) -- (m-3-4);

  \node[above right = -1pt of m-1-1] {\!\!(2,0)};
  \node[above right = -1pt of m-1-2] {\!\!(3,1)};
  \node[above right = -1pt of m-1-3] {\!\!(4,2)};
  \node[above right = -1pt of m-1-4] {\!\!(5,3)};
  
  \node[above right = -1pt of m-2-1] {\!\!(1,1)};
  \node[above right = -1pt of m-2-2] {\!\!(2,2)};
  \node[above right = -1pt of m-2-3] {\!\!(3,3)};
  \node[above right = -1pt of m-2-4] {\!\!(4,4)};

  \node[above right = -1pt of m-3-1] {\!\!(0,2)};
  \node[above right = -1pt of m-3-2] {\!\!(1,3)};
  \node[above right = -1pt of m-3-3] {\!\!(2,4)};
  \node[above right = -1pt of m-3-4] {\!\!(3,5)};

  
\end{tikzpicture}
\end{tabular}
\caption[LoF entry]{Minimal resolving sets (in red) for the complete graph $K_5$ (upper left), the path $P_5$ (upper right), the cycle $C_5$ (lower left), and the grid $G_{4, 3}$ (lower right).  Nodes are annotated with their distance vectors $d(u|R)$, which are all unique.}
\label{fig:path_and_complete_ex}
\end{figure}

For $G=(V,E)$ connected with $|V| = n \geq 2$, the path graph $P_n$, and the complete graph $K_n$ represent the two extremes of metric dimension. Indeed, $\beta(G) = 1$ if and only if $G = P_n$, and $\beta(G) = n-1$ if and only if $G = K_n$~\cite{chartrand2000resolvability}. Resolving sets for $P_n$ and $K_n$ are readily apparent (see \cref{fig:path_and_complete_ex}). 
For a path, either of the end vertices resolves every vertex, as each distance 0 to $(n-1)$ is attained exactly once. For a complete graph, every vertex is at distance 0 from itself and at distance 1 from all other vertices. This means that a single vertex $v \in V$ uniquely identifies itself but no other vertices. In order to distinguish all vertices, a resolving set of $K_n$ must be of size at least $(n-1)$. Moreover, any such set is resolving because the excluded vertex is the only one at a strictly positive distance from all other vertices.

The cycle graph $C_n$ on $n>2$ vertices has metric dimension 2~\cite{chartrand2000resolvabilityCycles}.
No single vertex set can resolve $C_n$ because every vertex has degree 2.
Next, we construct a set of size 2 and show that it is resolving.
Suppose that the vertex set of $C_n$ is $\{0,\ldots,n-1\}$, where consecutive integers are neighbors and so are $0$ and $n-1$. Let $R=\{0,1\}$. Then $d(x|R)=(\min\{x,n-x\},\min\{x-1,n+1-x\})$. But using that $\min\{a,b\}=(a+b-|a-b|)/2$, it follows that if $d(x|R)=d(y|R)$ then $|n-2x|=|n-2y|$ and $|n-2(x-1)|=|n-2(y-1)|$; in particular, $x$ and $y$ are at the same distance from $n/2$, and so are $x-1$ and $y-1$, which is possible only when $x=y$. Hence, $\beta(C_n)=2$ (\Cref{fig:path_and_complete_ex}).

Finally, consider the two-dimensional grid $G_{m,n}$ with dimensions $m,n \geq 1$.
The vertices of this graph correspond to ordered pairs $(x,y) \in \mathbb{Z}^2$ such that $0 \leq x < m$ and $0 \leq y < n$.
The edges correspond to pairs at Euclidean distance exactly 1, when considered as points in $\mathbb R^2$ (see~\cref{fig:path_and_complete_ex}).
In particular, every vertex has degree at least two and no singleton can be resolving.
Clearly, $d((0,0), (x,y)) = x+y$ and, more generally, $d((u,v), (x,y)) = |u-x| + |v-y|$.
Thus, shortest path distance in this case is equivalent to the $\ell_1$ norm or Manhattan distance.
Let $R = \{(0,0), (0,n-1)\}$. Observe that $d((a,b)|R) = (a+b, a+n-1-b)$. Then, for vertices $(a,b)$ and $(x,y)$, we have $d((a,b)|R) = d((x,y)|R)$ if and only if $a=x$ and $b=y$. Hence, $\beta(G_{m,n}) = 2$~\cite{khuller1996landmarks,melter1984metric}.
A symmetric argument shows that $R = \{(0,0), (m-1,0)\}$ is also resolving for $G_{m,n}$.

\section{Computational Complexity and Approximation}

Verifying that a given set of nodes $R$ in $G=(V,E)$ constitutes a resolving set is straightforward. For every $v \in V$, the vector of distances $d(v|R)$ can be generated in $O(|E|+|V| \log|V|)$ time. This collection of $\binom{|V|}{2}$ vectors then needs to be checked for duplicates. If all vectors are unique, the set is resolving, otherwise there is at least one pair of indistinguishable nodes in $G$ based on $R$.

A brute force solution determining the exact metric dimension of a general graph, on the other hand, requires an exhaustive search over a very large solution space. For a fixed set of size $s$, there are $\binom{|V|}{s}$ subsets of nodes that must be considered. Since we are interested in the smallest $s$ for which a subset of nodes of this size resolves $G$, increasing values of $s$ starting at 1 must be tested until a solution is found. Indeed, for a positive integer $k$, deciding whether $\beta(G) \leq k$ is an NP-complete problem. As a result, several approximation methods designed to find small resolving sets on general graphs have been developed. 

In what remains of this section we discuss one approach used to show the NP-completeness of the metric dimension problem, we describe an approximation algorithm based on a greedy selection criterion~\cite{hauptmann2012approximation}, and we outline two heuristic search techniques that have been applied to the problem of quickly finding small resolving sets.

\subsection{NP-Completeness}

The decision problem associated with metric dimension is to determine, given a graph $G$ and integer $k$, whether or not $\beta(G) \leq k$.
This decision problem is NP-complete, i.e., computationally intractable.
(For background on computational complexity, see, e.g., Goldreich~\cite{goldreich2008computational}.)
In this section, we present one proof of NP-completeness via reduction from 3-SAT, the problem of testing whether a given Boolean formula in conjunctive normal form, with three literals per clause, has a satisfying assignment~\cite{karp1972reducibility}. A reduction from the 3-dimensional matching problem to metric dimension is cited in~\cite{gary1979computers}, though we have not been able to find this proof in the literature.

Formally, the 3-SAT problem is as follows. Let $E$ be a Boolean expression in conjunctive normal form with $n$ variables $x_1, \dots, x_n$ and $m$ clauses $C_1, \dots, C_m$. For instance, the formula $(x_1 \vee \overline{x_2} \vee x_3) \wedge (x_2 \vee x_3 \vee \overline{x_4})$ consists of two clauses and four variables.
3-SAT is the problem of determining, given such a formula $E$, whether there exists an assignment mapping variables to truth values making $E$ true.
For the previous formula, setting $x_1$ and $x_2$ to True and $x_3$ and $x_4$ to False is one such an assignment.

\begin{figure}[h]
\centering
\tikzstyle{gadget}=[circle,thick,minimum size=4mm,draw=black]

\begin{tikzpicture}[every node/.style={gadget}]
  \node[label=$T_i$,fill=green!50] (T) at (0.6,0) {};
  \node[label=$a_i^1$] at (2,2) (a) {};
  \node[label=below:{$a_i^2$}] (b) at (2,1) {};
  \node[label=$b_i^1$] (c) at (4,2) {};
  \node[label=below:{$b_i^2$}] (d) at (4,1) {};
  \node[label=$F_i$,fill=blue!50] (F) at (5.4,0) {};

  \draw (T) -- (a) -- (c) -- (F) -- (d) -- (b) -- (T);
\end{tikzpicture}\qquad \tikzstyle{gadget}=[circle,thick,minimum size=4mm,draw=black]
\begin{tikzpicture}[every node/.style={gadget}]
  \node[label=$c_j^1$] (a) at (0,1) {};
  \node[label=$c_j^2$] (b) at (2,1) {};
  \node[label=$c_j^3$] (c) at (4,1) {};
  \node[label=left:$c_j^4$] (d) at (1,0) {};
  \node[label=right:$c_j^5$] (e) at (3,0) {};
  
  \draw (a) -- (b) -- (c);
  \draw (b) -- (d);
  \draw (b) -- (e);
\end{tikzpicture}
\caption[LoF entry]{Visualization of gadgets~\cite{khuller1996landmarks} for a single variable (left) and a single clause (right) in a Boolean expression.}
\label{fig:3sat_variable_and_clause}
\end{figure}
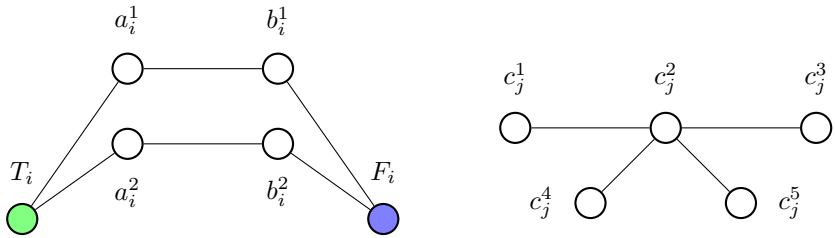

For an arbitrary 3-SAT instance $E$ we will construct a graph $G$ such that $E$ is satisfiable if and only if $\beta(G) = n+m$. We follow the construction of~\cite{khuller1996landmarks}. For every variable $x_i$ create a six cycle with nodes labeled $T_i$, $a^1_i$, $b^1_i$, $F_i$, $b^2_i$, and $a^2_i$, listed clockwise (see~\cref{fig:3sat_variable_and_clause} left). For every clause $C_j$ create a four star with nodes labeled $c^k_j$, with $1 \leq k \leq 5$, and central node $c^2_j$ (see~\cref{fig:3sat_variable_and_clause} right). 

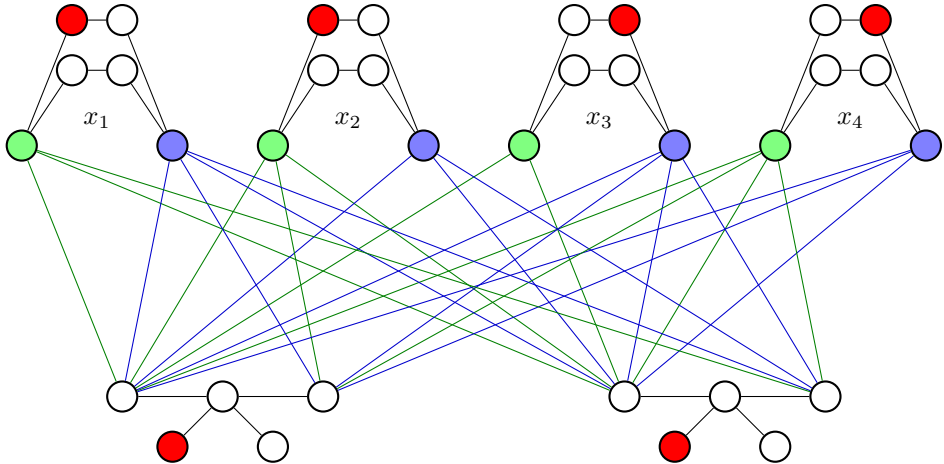
\begin{figure}[h]
\centering
\tikzstyle{gadget}=[circle,thick,minimum size=4mm,draw=black]
\begin{tikzpicture}[every node/.style={gadget},x=19pt,y=19pt]
  \node[fill=green!50] (T1) at (0,0) {};
  \node[fill=red] (a11) at (1,2.5) {};
  \node (a12) at (1,1.5) {};
  \node (b11) at (2,2.5) {};
  \node (b12) at (2,1.5) {};
  \node[fill=blue!50] (F1) at (3,0) {};
  \node[draw=white,minimum size=4mm] (x1) at (1.5,0.5) {$x_1$};
  \draw (T1) -- (a11) -- (b11) -- (F1) -- (b12) -- (a12) -- (T1);
  
  \def \s {5}
  \node[fill=green!50] (T2) at (\s+0,0) {};
  \node[fill=red] (a21) at (\s+1,2.5) {};
  \node (a22) at (\s+1,1.5) {};
  \node (b21) at (\s+2,2.5) {};
  \node (b22) at (\s+2,1.5) {};
  \node[fill=blue!50] (F2) at (\s+3,0) {};
  \node[draw=white,minimum size=4mm] (x1) at (\s+1.5,0.5) {$x_2$};
  \draw (T2) -- (a21) -- (b21) -- (F2) -- (b22) -- (a22) -- (T2);

  \def \s {10}
  \node[fill=green!50] (T3) at (\s+0,0) {};
  \node (a31) at (\s+1,2.5) {};
  \node (a32) at (\s+1,1.5) {};
  \node[fill=red] (b31) at (\s+2,2.5) {};
  \node (b32) at (\s+2,1.5) {};
  \node[fill=blue!50] (F3) at (\s+3,0) {};
  \node[draw=white,minimum size=4mm] (x1) at (\s+1.5,0.5) {$x_3$};
  \draw (T3) -- (a31) -- (b31) -- (F3) -- (b32) -- (a32) -- (T3);
  
  \def \s {15}
  \node[fill=green!50] (T4) at (\s+0,0) {};
  \node (a41) at (\s+1,2.5) {};
  \node (a42) at (\s+1,1.5) {};
  \node[fill=red] (b41) at (\s+2,2.5) {};
  \node (b42) at (\s+2,1.5) {};
  \node[fill=blue!50] (F4) at (\s+3,0) {};
  \node[draw=white,minimum size=4mm] (x1) at (\s+1.5,0.5) {$x_4$};
  \draw (T4) -- (a41) -- (b41) -- (F4) -- (b42) -- (a42) -- (T4);
  
  \def \d {-5}
  \node (c11) at (2,\d) {};
  \node (c12) at (4,\d) {};
  \node (c13) at (6,\d) {};
  \node[fill=red] (c14) at (3,\d-1) {};
  \node (c15) at (5,\d-1) {};
  \draw (c11) -- (c12) -- (c13);
  \draw (c12) -- (c14);
  \draw (c12) -- (c15);
  
  \def \s {10}
  \node (c21) at (\s+2,\d) {};
  \node (c22) at (\s+4,\d) {};
  \node (c23) at (\s+6,\d) {};
  \node[fill=red] (c24) at (\s+3,\d-1) {};
  \node (c25) at (\s+5,\d-1) {};
  \draw (c21) -- (c22) -- (c23);
  \draw (c22) -- (c24);
  \draw (c22) -- (c25);
  
  \foreach \s in {1,2}
  {
      \draw[green!50!black] (c\s1) -- (T1);
      \draw[green!50!black] (c\s1) -- (T2);
      \draw[green!50!black] (c\s1) -- (T3);
      \draw[green!50!black] (c\s1) -- (T4);
  }
  
  \draw[blue!80!black] (c11) -- (F1);
  \draw[blue!80!black] (c13) -- (F1);
  \draw[blue!80!black] (c11) -- (F3);
  \draw[blue!80!black] (c13) -- (F3);
  
  \draw[blue!80!black] (c21) -- (F2);
  \draw[blue!80!black] (c23) -- (F2);
  \draw[blue!80!black] (c21) -- (F3);
  \draw[blue!80!black] (c23) -- (F3);
  
  \draw[blue!80!black] (c11) -- (F2);
  \draw[green!50!black] (c13) -- (T2);
  \draw[blue!80!black] (c21) -- (F4);
  \draw[green!50!black] (c23) -- (T4);
  
  \draw[blue!80!black] (c11) -- (F4);
  \draw[blue!80!black] (c13) -- (F4);
  \draw[green!50!black] (c13) -- (T4);
  
  \draw[blue!80!black] (c21) -- (F1);
  \draw[blue!80!black] (c23) -- (F1);
  \draw[green!50!black] (c23) -- (T1);
\end{tikzpicture}
\caption[LoF entry]{The graph produced by the reduction for the formula $E = (x_1 \vee \overline{x_2} \vee x_3) \wedge (x_2 \vee x_3 \vee \overline{x_4})$.  Gadgets for variables and clauses are oriented as in Figure~\ref{fig:3sat_variable_and_clause}.  As $E$ is satisfiable, the graph has a resolving set of size $4+2=6$, shown in red.}
\label{fig:3sat_clause_ex}
\end{figure}

These cycles and stars are used to form a connected graph by including the edge $\{T_i,c^1_j\}$ for every variable $x_i$ and every clause $C_j$. In addition, when $x_i$ is used as a positive literal in $C_j$, the edges $\{F_i,c^1_j\}$ and $\{F_i,c^3_j\}$ are added to the graph. Instead, when $x_i$ is used as a negative literal in $C_j$ the edges $\{F_i,c^1_j\}$ and $\{T_i,c^3_j\}$ are added. Otherwise, if $x_i$ does not appear in $C_j$ the edges $\{F_i,c^1_j\}$, $\{F_i,c^3_j\}$, and $\{T_i,c^3_j\}$ are added. (See~\cref{fig:3sat_clause_ex}.) Since there are a total of $6n+5m$ nodes in the final graph, and there are at most four edges between the subgraphs representing $x_i$ and $C_j$ for all $1 \leq i \leq n$ and $1 \leq j \leq m$, this construction takes polynomial time.

Notice that any resolving set $R$ of $G$ must include at least one of $\{a^1_i,a^2_i,b^1_i,b^2_i\}$ for all $1 \leq i \leq n$ and at least one of $\{c^4_j,c^5_j\}$ for all $1 \leq j \leq m$. Hence $\beta(G) \geq n+m$.
Showing that $\beta(G) = n+m$ if $E$ is satisfiable is straightforward. 
Given an assignment of the variables to true or false, take $R = \{c^4_j \mid 1\leq j\leq m\}\cup\{a^1_i\mid x_i\text{ is true}\}\cup\{b^1_i\mid x_i\text{ is false}\}$. (One could also choose $c^5_j$, $a^2_i$, and/or $b^2_i$).

To see that $R$ is resolving,
first consider the variable gadget for $x_i$, and any clause $C_j$.
Recall that $c^4_j\in R$.
The vertices of the $x_i$ gadget are split into two groups based on distances to $c^4_j$: vertices $\{T_i,F_i\}$ at distance 3 and $\{a^1_i,a^2_i,b^1_i,b^2_i\}$ at distance 4.
Any vertex in the set $\{a^1_i,a^2_i,b^1_i,b^2_i\}$, therefore, serves to disambiguate the elements of these groups: $T_i$ and $F_i$ will have distance 1 and 2, or vice versa, and the remaining nodes will attain every distance in $\{0,1,2,3\}$.
As we have either $a^1_i\in R$ or $b^1_i\in R$, the variable gadget is resolved.
In fact, this statement holds regardless of whether the formula is satisfiable.

Now consider a clause $C_j$ and a variable $x_i$ that causes this clause to be satisfied. 
There are two cases: either $x_i$ is a positive literal in $C_j$ and is given a value of true, or $x_i$ is a negative literal in $C_j$ and is given a value of false.
In the first case, we have $a^1_i \in R$, and recall $c^4_j\in R$.
In the $C_j$ gadget, $c^2_j$ is the unique vertex at distance 1 from $c^4_j$, and $c^5_j$ is the unique vertex at distance 4 from $a^1_i$, distinguishing these vertices from all others. 
Finally, we have $d(a^1_i, c^1_j) = 2$ and $d(a^1_i, c^3_j) = 3$, finishing the proof of the first case.
The second case is symmetric to the first, with the edge $\{F_i,c^3_j\}$ becoming $\{T_i,c^3_j\}$, and the role of $a^1_i$ played by $b^1_i$ instead.

Conversely, it can be shown that $E$ is satisfiable if $\beta(G) = n+m$ by setting $x_i$ to true if either $a^1_i$ or $a^2_i$ is in the resolving set and to false otherwise.
Thus, this construction reduces 3-SAT to the metric dimension decision problem in polynomial time.

Although metric dimension is a computationally difficult problem on arbitrary graphs, there are efficient and fixed-parameter tractable algorithms in certain restricted settings. Indeed, trees~\cite{harary1976metric,slater1975leaves}, cographs~\cite{epstein2015weighted}, and outerplanar graphs~\cite{diaz2012complexity} admit linear or polynomial time algorithms, and interval~\cite{foucaud2017identification2} and permutation graphs~\cite{belmonte2015metric} have fixed-parameter tractable algorithms with respect to resolving set size. There are also fixed-parameter tractable algorithms for metric dimension on general graphs with respect to other graph quantities including vertex cover~\cite{hartung2013parameterized}, maximum leaf number~\cite{eppstein2015metric}, and modular width~\cite{belmonte2015metric}.

\subsection{Information Content Heuristic}

Approximate solutions to the metric dimension problem on general graphs may be generated using a fairly simple greedy algorithm. Originally devised for use on instances of the test set problem~\cite{berman2005tight} and later modified for estimating the metric dimension of graphs~\cite{hauptmann2012approximation}, this algorithm is based on an Information Content Heuristic (ICH).

The underlying concept is similar to information gain as it relates to the construction of decision trees~\cite{quinlan1986induction}.
Let $G=(V,E)$ be a graph with $n = |V|$.
Consider the classification problem for which each vertex is a training example from a unique class. The feature vector of each $v \in V$ is simply $d(v|V) \in \mathbb{R}^n$, the shortest path distances to all vertices in the graph.
In particular, the $(n \times n)$ distance matrix $\mathbf{D}$ associated with $G$ fully describes the training data.
The choice of resolving set $R$ can thus be thought of as the choice of some subset of features: $d(v|R)$ is simply a selection of entries of $d(v|V)$ corresponding to the ``features'' $R$.

In the usual greedy algorithm to construct a decision tree, at each node of the decision tree, one selects the feature which maximizes the information gained about the true class at the child nodes.
In the case of resolving sets, each node of the decision tree corresponds to the addition of some vertex $u$ to $R$, with branches corresponding to the possible values of $d(\cdot,u)$.
To measure information in our setting, consider the distribution $p_R$ induced by the equivalence classes of $d(\cdot|R)$, which assigns probability $|\{v \in V : d(v|R) = d\}|/n$ to each possible value $d \in \{d(v|R) : v\in V\}$.
Then the information of $R$ is measured by the Shannon entropy of $p_R$, denoted $H(R)$.
Letting $R_t$ be the resolving set at iteration $t$, with $R_0 = \emptyset$, the ICH algorithm therefore chooses the vertex $v_t\in V$ to maximize $H(R_t\cup\{v_t\})$.
The algorithm terminates when $H(R_t) = \log n$, the maximum possible entropy over $n$ items, indicating that all $n$ vertices are uniquely represented by their distances to $R_t$.\footnote{In fact, the ICH algorithm is exactly the same as the information gain algorithm for decision trees, under the constraint that the decision made at every node of the tree at the same level must be the same, i.e., we choose the same vertex at iteration $t$ for all decision nodes at level $t-1$.
To see the equivalence, let $n_d = |\{v \in V : d(v|R) = d\}|$ be the size of equivalence class $d$, and $p_d$ be the distribution of labels within equivalence class $d$, i.e., the unform distribution on all $n_d$ vertices in the equivalence class.
Then maximizing information gain, weighted by the size of each decision node (recall that we must choose the same feature to split at all nodes), is the same as maximizing $H(p_R)$:
$\min \sum_{d} \frac{1}{n_d} H(p_d) = \min \sum_{d} \frac{1}{n_d} \log(n_d) = \max - \sum_{d} \frac{1}{n_d} \log(\frac{1}{n_d}) = \max H(p_d)$.}

Asymptotically, the time complexity of the ICH algorithm is $O(n^3)$. This makes it an effective algorithm only for approximating the metric dimension of relatively small graphs. Nevertheless, this algorithm does have the significant advantage of guaranteeing a $1 + (1 + o(1))\cdot\ln(n)$ approximation ratio (i.e., the approximate metric dimension of a graph discovered by the ICH is never more than $1 + (1 + o(1))\cdot\ln(n)$ times as large as its true metric dimension). This is the best possible approximation ratio for the metric dimension problem~\cite{hauptmann2012approximation}.

\subsection{Other Heuristics}

Among many heuristic methods commonly deployed in non-convex search problems, genetic algorithms and variable neighborhood search in particular have been used to find small resolving sets on general graphs with some success.
Genetic algorithms, inspired by the concept of biological evolution, seek optimal solutions to problems by incrementally changing a population of candidate solutions from one generation to the next through the biologically motivated operations of mutation and selective recombination~\cite{davis1991handbook}. This approach has been shown to perform quite well when applied to metric dimension in comparison to other state-of-the-art algorithms, including methods based on an integer programming formulation of the problem and the CPLEX~\cite{cplex20059} optimization package~\cite{kratica2009computing}.

The variable neighborhood search (VNS) technique starts with an initial, non-optimal solution and iteratively expands a neighborhood on which to perform a local search. When a point in the space which improves upon the initial solution is found, the search is restarted with this point at its center. In the context of searching for small resolving sets, VNS seems to outperform genetic algorithm based methods on many kinds of graphs and has been used to improve upon previous upper bounds for certain hypercubes~\cite{mladenovic2012variable}.

\section{Graph Features and Metric Dimension Relationships}

In this section, we overview some general observations about metric dimension and its relationship to other graph quantities.
These observations are often useful in bounding or exactly determining the metric dimension of a given graph.

\subsection{Diameter and Metric Dimension} The diameter of a graph $G=(V,E)$, denoted $\diam(G)$, is the length of a longest shortest path in $G$. For ease of notation, let $\delta=\diam(G)$ and $k=\beta(G)$. It is not surprising that a relationship exists between $\delta$ and $k$. Indeed, let $R \subset V$ be a minimum resolving set of $G$ and consider $d(v|R)$ for each $v \in V$. Since such vectors can only contain a 0 when $v \in R$ and, for $v\notin R$, $1\le d(r,v)\le\delta$, it follows that $|V| \leq \delta^k+k$~\cite{khuller1996landmarks}.
This bound is usually loose, though graphs with $|V|=\delta+k$ have been fully characterized~\cite{hernando2010extremal}. 
The related bound,
\begin{displaymath}n \leq (\lfloor 2\delta/3 \rfloor + 1)^k + k \sum_{i=1}^{\lceil \delta/3 \rceil} (2i-1)^{k-1},\end{displaymath} 
is generally tighter~\cite{hernando2010extremal}.
This bound can be made tighter still for specific families of graphs. We list a few results here~\cite{beaudou2018bounding}.
\begin{itemize}
    \item If $G$ is a tree, $|V| \leq (\delta k+4)(\delta+2)/8$, with equality for trees with even diameter~\cite{beaudou2018bounding}.
    \item If $G$ is an outerplanar graph, $|V| = O(\delta^2k)$~\cite{beaudou2018bounding}.
    \item If $K_i$ is not a minor of $G$, $|V| \leq (\delta k+1)^{i-1}+1$~\cite{beaudou2018bounding}.
    \item If $G$ has constant treewidth, $|V| = O(k\delta^{O(1)})$~\cite{beaudou2018bounding}.
    \item If the rankwidth of $G$ is at most $r$, $|V| \leq (\delta k+1)^{\delta(3(2^r)+2)}+1$~\cite{beaudou2018bounding}.
    \item If $G$ is an interval or permutation graph, $|V| = O(\delta k^2)$~\cite{foucaud2017identification}.
    \item If $G$ is a unit interval graph, $|V| = O(\delta k)$~\cite{foucaud2017identification}.
\end{itemize}

\subsection{Twin Nodes and Metric Dimension}
\label{sec:twin_nodes}
Let $G=(V,E)$ be an undirected graph and, for each $v\in V$, define the closed-neighborhood of $v$ as $N(v) = \{u \mid \{u,v\} \in E\}$. We call $u,v \in V$ twins when $N(u) \cup \{u\} = N(v) \cup \{v\}$.

Twin nodes have an interesting relationship to metric dimension. In fact, when $u$ and $v$ are twins, $d(u,w) = d(v,w)$, for all $w \in V \setminus \{u,v\}$. As a result, any resolving set of $G$, minimal or not, must include at least one of $u$ and $v$. More formally, define over $V$ the equivalence relation: $u \equiv v$ if and only if $u$ and $v$ are twins. Let $\tau(G)$ be the set of twin equivalence classes of $G$. Then, if $R$ is a resolving set of $G$, $|R \cap \tau| \geq |\tau|-1$, for each $\tau \in \tau(G)$. In particular~\cite{hernando2010extremal}:
\begin{displaymath}\beta(G) \geq \sum_{\tau \in \tau(G)} (|\tau|-1).\end{displaymath}
Twin nodes have been used, for example, to study connections between metric dimension, diameter, and graph size~\cite{hernando2010extremal}.

\subsection{Graphs with Extreme Metric Dimension} As we saw in~\Cref{sec:simple_examples}, for $G=(V,E)$, $1 \leq \beta(G) \leq n-1$, where $\beta(G) = 1$ if and only if $G \cong P_n$, and $\beta(G)=n-1$ if and only if $G \cong K_n$. In fact, for every $1 \leq k \leq n-1$ there is a connected graph $G$ with $n$ vertices and $\beta(G) = k$~\cite{chartrand2000resolvability}. Suppose $G=(V,E)$ such that $\beta(G) = 2$. Such graphs have not been fully characterized but there are a set of simple properties that they must have~\cite{khuller1996landmarks}. In particular, if $\{u,v\} \subset V$ is a resolving set of $G$ of minimum size:
\begin{enumerate}
    \item $G$ cannot contain $K_5$ as a subgraph.
    \item $G$ cannot contain $K_{3,3}$ as a subgraph. ($K_{m,n}$ denotes the complete bipartite graph with partitions of size $m$ and $n$.)
    \item There is a unique shortest path between $u$ and $v$.
    \item $\deg(w) \leq 5$ for all nodes $w \in V$ on the shortest path between $u$ and $v$.
    \item $\deg(u) \leq 3$ and $\deg(v) \leq 3$.
\end{enumerate}

\noindent Properties $(1)$ and $(2)$ may bring Wagner's theorem~\cite{wagner1937eigenschaft} to mind, a characterization of planar graphs forbidding $K_5$ and $K_{3,3}$ as minors, suggesting that $G$ must be planar.
On the contrary, there are non-planar graphs with metric dimension 2~\cite{khuller1996landmarks}.

Graphs $G=(V,E)$ with $\beta(G) = (n-2)$, on the other hand, have been fully characterized~\cite{chartrand2000resolvability}. For two graphs $G_1$ and $G_2$, let $G_1 \cup G_2$ denote their disjoint union, and let $G_1 + G_2$ denote the graph formed by taking a disjoint union and joining every node in $G_1$ with every node in $G_2$. Furthermore, define $\overline{K}_n$ to be a graph with $n$ nodes and no edges. Then the metric dimension of a graph with $n$ nodes is $n-2$ if and only if the graph is one of the following:

\begin{itemize}
    \item A complete bipartite graph, $K_{s,t}$ with $s,t \geq 1$.
    \item $K_s + \overline{K}_t$ with $s \geq 1$ and $t \geq 2$.
    \item $K_s + (K_1 \cup K_t)$ with $s, t \geq 1$.
\end{itemize}

These characterizations of graphs with large and small metric dimension, along with relationships to other graph quantities like total twin nodes and diameter, are often useful when beginning the search for the metric dimension of specific kinds of graphs.

\section{Specific Families of Graphs}

While determining the metric dimension of arbitrary graphs is a computationally complex task, exact formulae, upper bounds, and polynomial time algorithms exist for certain types of graphs.
In practical applications these analytic and algorithmic results are of crucial importance.
Generic approximation algorithms tend to provide small resolving sets but do not scale well enough with network size to be useful on larger networks which routinely include more than $10^6$ total nodes~\cite{cha2010measuring,meusel2014graph,TilLla18}.
Tailored algorithms for specific graph structures can therefore be useful in quickly discovering small resolving sets even on large networks. 

In what follows, we discuss some of the more prominent graph families on which metric dimension has been studied. \Crefrange{tab:beta_list_exact}{tab:beta_list_other} provide a short, incomplete list of known bounds with references.
The literature also addresses Cartesian products of graphs~\cite{caceres2007metric,jiang2019metric} and infinite graphs~\cite{caceres2009metric}, which we do not discuss.

\subsection{Fans and Wheels}
\label{sec:fans-wheels}
Recall the path $P_n$ and cycle $C_n$ from \Cref{sec:simple_examples}, two of the simplest graphs to study metric dimension.
Fans and wheels are simple modifications of these which add a new fully-connected vertex.
Formally, the fan graph on $n+1$ vertices, denoted $F_n$, consists of a path $P_n$ on $n$ vertices and one additional vertex, $a$, adjacent to all vertices on $P_n$.
Similarly, the wheel graph $W_n$ is a cycle $C_n$ of size $n$ with an additional vertex, $a$, adjacent to all vertices on $C_n$.
As the metric dimensions of graphs $P_n$ and $C_n$ are elementary to determine (1 and 2, respectively), one might expect the metric dimensions of $F_n$ and $W_n$ to be similarly trivial.
Yet despite their simplicity, fans and wheels have much more complex expressions for their metric dimension.
Part of this complexity stems from the diameter of these graphs shrinking from order $n$ to at most 2, reducing the possible shortest path distances to the range $\{0,1,2\}$. Specifically, their metric dimensions are given as follows~\cite{shanmukha2002metric}.
\begin{equation}
\beta(F_{x+5k}) = 
\beta(W_{x+5k}) =
\begin{cases}
3+2k, & \text{ when } x \in \{7,8\}, \\
4+2k, & \text{ when } x \in \{9,10,11\},
\end{cases}
\label{ide:betaPFnWn}
\end{equation}
for all $k \geq 0$. We find instructive to provide a high level proof of this result. Accordingly, we assume in what follows that $n>6$.

Focusing first on $F_n$, observe that $a$ appears in none of its minimal resolving sets.
To see this, suppose for a contradiction that $R$ is a minimal resolving set of $F_n$ such that $a \in R$.
Since $R$ is minimal, there must be two vertices $u,v\in V$ with $d(u|R\setminus\{a\}) = d(v|R\setminus\{a\})$. 
Furthermore, without loss of generality we must have $u=a$, since otherwise $d(u,a)=d(v,a)=1$ and $R$ would not resolve $u$ and $v$. 
Thus, as $u=a$, and as $d(a,R\setminus\{a\})$ is the all-ones vector, $v$ must be distance 1 from every node in $R\setminus\{a\}$.
By definition of $F_n$, there are at most two nodes other than $a$ at distance 1 from $v$, so $|R\setminus\{a\}| \leq 2$. 
One can check that there are now only 6 possible values of $d(\cdot|R\setminus\{a\})$, namely $(1,2), (0,2), (1,1), (2,0), (2,0), (2,1), (2,2)$.
As $a$ is at distance 1 from all nodes on the path, and $n>6$ we cannot have resolved $F_n$. 
Hence, no minimal resolving set of $F_n$ contains the vertex $a$.

Let $R$ be a minimal resolving set of $F_n$ of size $\beta$. Note that at most one vertex of $P_n$ can be at distance 2 from all vertices in $R$. Besides, at most $\beta$ vertices of $P_n$ can be adjacent to exactly one vertex of $R$. 
Otherwise, if there were more, the pigeonhole principle would imply that that at least two vertices are adjacent to the same vertex of $R$ and, therefore, would be indistinguishable. As a result, all the remaining nodes of $P_n$ must either be in $R$, or be adjacent to exactly two vertices of $R$. Note that having fewer than $\beta$ vertices adjacent to a single element of $R$ cannot increase the overall number of vertices resolved. 
Hence $n\le 1+2\beta+\beta/2$, or equivalently: $\lceil2(n-1)/5\rceil\le\beta$. But, because $n>6$, we may write $n=x+5k$ with $7\le x\le 11$ and $k\ge0$. In particular~\cite{shanmukha2002metric}:
\begin{displaymath}
\beta(F_{x+5k})\geq
\begin{cases}
3+2k, & \text{ when } x \in \{7,8\}; \\
4+2k, & \text{ when } x \in \{9,10,11\}.
\end{cases}
\end{displaymath}

\begin{figure}[h]
  \centering
\tikzstyle{point}=[circle,minimum size=1pt, inner sep=0.4pt,draw=black, fill=black]
\tikzstyle{group}=[draw=black,dashed,rounded corners,fill=blue,fill opacity=0.05]
\begin{tikzpicture}[scale=0.9]
    \node[non-res] (0) at (6,1.8) {};
    \node[non-res] (1) at (0,0) {};
    \node[res] (2) at (1,0) {};
    \node[non-res] (3) at (2,0) {};
    \node[res] (4) at (3,0) {};
    \node[non-res] (5) at (4,0) {};
    \node[non-res] (6) at (5,0) {};
    \node[res] (7) at (6,0) {};
    \node[non-res] (8) at (7,0) {};
    \node[res] (9) at (8,0) {};
    \node[non-res] (10) at (9,0) {};
    \node[point] (a) at (9.5,0) {};
    \node[point] (b) at (9.75,0) {};
    \node[point] (c) at (10,0) {};
    \node[non-res] (11) at (10.5,0) {};
    \node[non-res] (12) at (11.5,0) {};
    \node[res] (13) at (12.5,0) {};
    \draw (1) -- (2) -- (3) -- (4) -- (5) -- (6) -- (7) -- (8) -- (9) -- (10);
    \draw (11) -- (12) -- (13);

    \foreach \i in {1,...,13}
    {
        \draw (0) -- (\i);
    }

    \node[group,fit=(1)(5)] {};
    \node[group,fit=(6)(10)] {};
    \node[group,fit=(11)(13)] {};
\end{tikzpicture} 
  \caption[LoF entry]{A visualization of a minimal resolving set of the fan $F_n$ when $(n \bmod 5) = 3$. Dashed lines group separate consecutive blocks of five nodes.}
  \label{fig:fan_min_ex}
\end{figure}
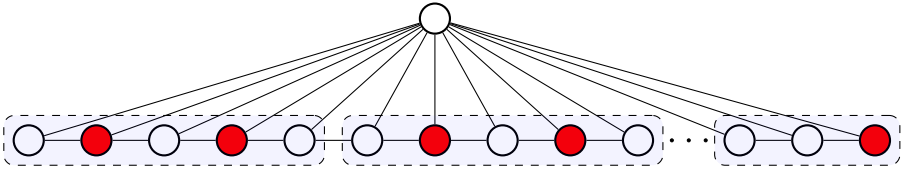

Furthermore, it is easy to see that this lower bound on $\beta(F_{x+5k})$ is also an upper bound via a simple construction. In particular, there is a resolving set of $F_{x+5k}$ of size $3+2k$ when $x$ is 7 or 8 and of size $4+2k$ when $x$ is 9, 10, or 11. Let $m = n \bmod 5$ and let $R$ be a minimal resolving set of $F_n$. For each consecutive, full block of five vertices, $v_x, v_y \in R$ where $x = 5j+2$ so that $v_x$ is the second vertex of the $j^{th}$ block and $y = 5j+4$ so that $v_y$ is the fourth vertex of the $j^{th}$ block for $j\geq 0$. If $m \in \{2,3\}$, $v_n \in R$. If instead $m = 4$, $v_{(n-2)}, v_n \in R$ (see Figure~\ref{fig:fan_min_ex}).

Finally, we claim that $\beta(W_n) = \beta(F_n)$ for $n>6$. Indeed, suppose without loss of generality that $R$ is a minimal resolving set of $F_n$ such that $1,n \not\in R$. (Any resolving set of $F_n$ can be made to satisfy this requirement by replacing $1$ with $2$ or $3$, and $n$ with $n-1$ or $n-2$, depending on $R \setminus \{1,n\}$.) Since the only difference between the fan and wheel graph is the inclusion of the edge $\{1,n\}$ in $W_n$, and since this edge is not required to determine shortest path distances between elements of $R$ and any other vertices in $W_n$, $R$ is also a resolving set for $W_n$. So $\beta(W_n) \leq \beta(F_n)$. Conversely, suppose that $R$ is a minimal resolving set of $W_n$; in particular, $|R|\leq\beta(F_n)$. Then there must be at least one edge $\{i,j\}$ in $W_n$ such that $i,j \not \in R$ and $i,j \neq a$. This follows from $\beta(W_n) < \lceil \frac{n-1}{2} \rceil$. 
Removing this edge, therefore, does not affect $R$ as a resolving set. In particular, $R$ also resolves $F_n$ and $\beta(W_n) = \beta(F_n)$, which shows~(\ref{ide:betaPFnWn}).

\subsection{Trees and Unicyclic Graphs}
\label{sec:trees-unicyclic-graphs}
Problems that are computationally difficult on general graph structures often admit more efficient solutions on trees. This is the case for metric dimension too. A simple formula giving the metric dimension of trees that are not also paths~\cite{chartrand2000resolvability,harary1976metric,slater1975leaves} leads immediately to a polynomial time algorithm for finding resolving sets of minimum size on trees. To begin, we present several important definitions. Let $G=(V,E)$ be a general graph, not necessarily a tree. In what follows, for $v\in V$, we use $\deg(v) = |\{u \in V\mid (v,u)\in E\}|$ to denote the degree of $v$.

\begin{definition}{(Leaf Vertex.)}
A vertex $\ell \in V$ is called a leaf when $\deg(\ell) = 1$. The number of leaves in $G$ is denoted $\ell(G)$.
\end{definition}

\begin{definition}{(Major Vertex and Terminal Degree.)}
A vertex $v \in V$ is called a major vertex when $\deg(v) \geq 3$.
The terminal degree of a major vertex $v \in V$ is the number of leaves $\ell \in V$ such that $d(\ell, v) < d(\ell, u)$, for all other major vertices $u \in V$.
\end{definition}

\begin{definition}{(Exterior Major Vertex.)}
A major vertex of $G$ is called exterior when its terminal degree is strictly positive. The number of exterior major vertices in $G$ is denoted $ex(G)$.
\end{definition}

\begin{figure}[h]
\centering
\begin{tikzpicture}[scale=0.8,level/.style={sibling distance = 3cm/#1,
  level distance = 1.25cm},every node/.append style={font=\scriptsize}] 
  \node[non-res] {1}
    child { node [non-res] {2} 
        child { node [non-res] {6} 
            child { node [non-res] {12} }}
        child { node [non-res] {7} 
            child { node [res] {13}}
            child { node [non-res] {14}}}}
    child { node [non-res] {3} 
        child { node [res] {8}}}
    child { node [non-res] {4}}
    child { node [non-res] {5} 
        child { node [non-res] {9}
            child { node [res] {15}}
            child { node [res] {16}}
            child { node [non-res] {17}}}
        child { node [non-res] {10}}
        child { node [non-res] {11}
            child { node [res] {18}}}};
\end{tikzpicture}
\caption[LoF entry]{A tree of size 18. Vertices 4, 8, 10, 12, 13, 14, 15, 16, 17, and 18 are leaves. The vertices 1, 2, 5, 7, and 9 are exterior major vertices with terminal degree 2, 1, 2, 2, and 3, respectively. $R = \{8, 13, 15, 16, 18\}$ is a resolving set of minimum size.}
\label{fig:tree}
\end{figure}
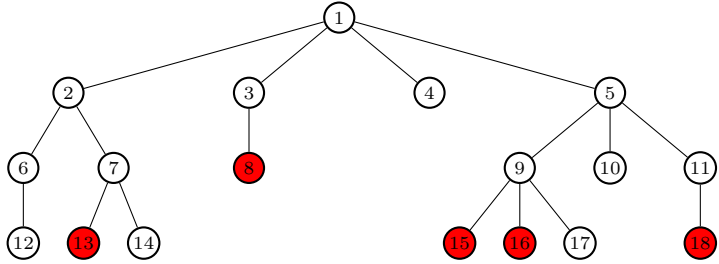

From these definitions, we can write $\beta(G) \geq \ell(G) - ex(G)$, with equality when $G$ is a tree such that $ex(G) > 0$, or equivalently, when $G$ is not a path~\cite{chartrand2000resolvability}. 
Moreover, any set $R \subset V$ which contains every leaf, except one, associated with each exterior major vertex is a subset of a minimal resolving set in $G$. When $G$ is a tree, any such $R$ is resolving~\cite{chartrand2000resolvability}.
These observations permit an $O(|V| + |E|)$ algorithm for constructing minimal resolving sets on trees: after partitioning the leaves of a tree based on exterior major vertices using a depth first search, one element of each partition may be dropped to produce a resolving set of minimum size (see~\Cref{fig:tree}).

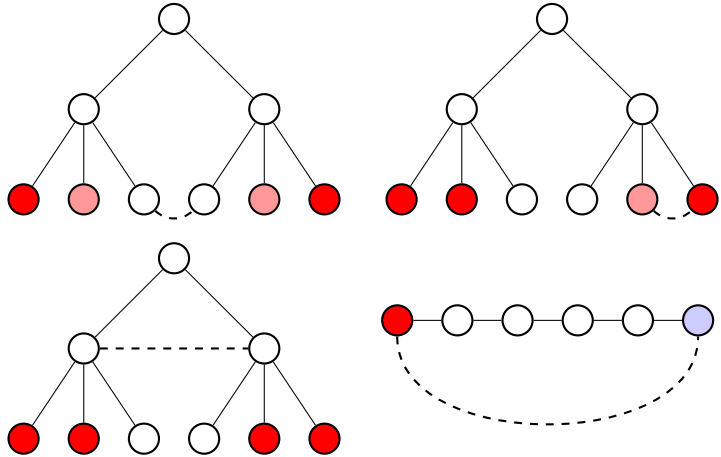
\begin{figure}[h]
  \tikzstyle{new-edge}=[dashed, thick]
  \begin{tabular}{cc}
    \begin{tikzpicture}[scale=0.8,level 1/.style={sibling distance=30mm},level 2/.style={sibling distance=10mm}]
  \node[non-res] (1) {}
    child { node [non-res] (2) {}
        child { node [res] (4) {} }
        child { node [light-res] (5) {} }
        child { node [non-res] (6) {} }}
    child { node [non-res] (3) {}
        child { node [non-res] (7) {} }
        child { node [light-res] (8) {} }
        child { node [res] (9) {} }};
    \draw[new-edge] (6) to [out=315,in=225] (7);
\end{tikzpicture} & \begin{tikzpicture}[scale=0.8,level 1/.style={sibling distance=30mm},level 2/.style={sibling distance=10mm}]
  \node[non-res] (1) {}
    child { node [non-res] (2) {}
        child { node [res] (4) {} }
        child { node [res] (5) {} }
        child { node [non-res] (6) {} }}
    child { node [non-res] (3) {}
        child { node [non-res] (7) {} }
        child { node [light-res] (8) {} }
        child { node [res] (9) {} }};
    \draw[new-edge] (8) to [out=315,in=225] (9);
\end{tikzpicture} \\[4pt]
    \begin{tikzpicture}[scale=0.8,level 1/.style={sibling distance=30mm},level 2/.style={sibling distance=10mm}]
  \node[non-res] (1) {}
    child { node [non-res] (2) {}
        child { node [res] (4) {} }
        child { node [res] (5) {} }
        child { node [non-res] (6) {} }}
    child { node [non-res] (3) {}
        child { node [non-res] (7) {} }
        child { node [res] (8) {} }
        child { node [res] (9) {} }};
    \draw[new-edge] (2) to [out=0,in=180] (3);
\end{tikzpicture} & \begin{tikzpicture}[scale=0.8]
    \node[res] (1) at (0,0) {};
    \node[non-res] (2) at (1,0) {};
    \node[non-res] (3) at (2,0) {};
    \node[non-res] (4) at (3,0) {};
    \node[non-res] (5) at (4,0) {};
    \node[blue] (6) at (5,0) {};
    \draw (1) -- (2) -- (3) -- (4) -- (5) -- (6);
    \draw[new-edge] (1) to [out=270,in=270] (6);
\end{tikzpicture}
  \end{tabular}
  \centering
\caption[LoF entry]{Example of trees $T$ such that, when a single edge (dashed) is added to form the unicyclic graph $G$, $\beta(G)=\beta(T)-2$ (top left), $\beta(G) = \beta(T)-1$ (top right), $\beta(G) = \beta(T)$ (bottom left), or $\beta(G) = \beta(T)+1$ (bottom right)~\cite{chartrand2000resolvability}. In each example, red nodes belong to a resolving set for both $G$ and $T$, light red nodes are used to resolve $T$, and light blue nodes are used to resolve $G$.}
\label{fig:unicyclic_examples}
\end{figure}

Let $G=(V,E)$ be a unicyclic graph (i.e. a graph that can be expressed as a tree with a single additional edge) with $|V|\geq3$. 
Let $T$ be any spanning tree of $G$ and $e$ the only edge in $G$ that is not in $T$. Then $\beta(T)-2 \leq \beta(G) \leq \beta(T)+1$~\cite{chartrand2000resolvability,poisson2002dimension}. As illustrated in~\Cref{fig:unicyclic_examples}, there are unicyclic graphs achieving each of the values in the integer interval $[\beta(T)-2, \beta(T)+1]$.

To justify the lower bound on $\beta(G)$ for unicyclic $G$, we consider three cases.
First, if $e$ is incident on leaves in $T$ then $\ell(G)=\ell(T)-2$, and $ex(G) \leq ex(T)$ because $G$ and $T$ have the same major vertices but $e$ reduces the terminal degree of at least one major vertex in $T$.
In particular, $\beta(G) \geq \ell(G) - ex(G) \geq (\ell(T) - 2) - ex(T) = \beta(T) - 2$. Instead, if $e$ is incident on exactly one leaf in $T$ then $\ell(G)=\ell(T)-1$, and $ex(G) \leq ex(T) + 1$ because $e$ may turn a vertex in $T$ into an exterior major vertex. So, $\beta(G) \geq \ell(G) - ex(G) \geq (\ell(T) - 1) - (ex(T)-1) = \beta(T) - 2$. Finally, if $e$ is not incident on any leaf in $T$, $ex(G) \leq ex(T)+2$ because both vertices $e$ is incident on may become exterior major vertices. Hence, $\beta(G) \geq \ell(G) - ex(G) \geq \ell(T) - (ex(T)+2) = \beta(T) - 2$.

Slightly more work is required to verify the upper bound on the metric dimension of unicyclic graphs. One approach is to focus on a subset of major vertices on the cycle in $G$. Let $W$ contain every leaf in $G$, except one, associated with each exterior major vertex, and let $m$ be the number of major vertices on the cycle in $G$ with a branch to an element of $W$. Proceeding by cases, it can be shown that $\beta(G) \leq \beta(T)+1$ whether $m\geq 3$, $m=2$, $m=1$, or $m=0$~\cite{chartrand2000resolvability}.

\subsection{Grids, Honeycombs, and Hexagon Networks} In~\cref{sec:simple_examples}, we examined two-dimensional grids, $G_{m,n}$, and argued that $\beta(G_{m,n}) = 2$.
For the $d$-dimensional grid $G_{n_1,\dots,n_d}$ defined in an analogous manner for $n_i>1$ for $1\leq i \leq d$, $\beta(G_{n_1, \dots, n_d}) = d$~\cite{khuller1996landmarks}.
The key idea is that if $r_0 := (0, \dots, 0)$, and $r_i$ is the vector of zeroes with entry $(n_i-1)$ in the $i$-th position, then for any vertex $v = (x_1, \dots, x_d)$ in the grid we have:
\begin{align*}
    d(r_0, v) &= x_1 + \dots + x_d \\
    d(r_1, v) &= (n_1 - 1 - x_1) + x_2 + \dots + x_d \\
    \cdots \\
    d(r_{d-1}, v) &= x_1 + \dots + (n_{d-1} - 1 - x_{d-1}) + x_d.
\end{align*}
Since this linear system of $d$ equations and $d$ unknowns is invertible, $\{r_0, \dots, r_{d-1}\}$ resolves the grid i.e. $\beta(G_{n_1,\dots,n_d})\le d$. The claim follows after showing that the linear system associated with any set $R$ with fewer than $d$ elements is not invertible.

\begin{figure}[h]
\centering
\tikzstyle{point}=[circle,minimum size=1pt, inner sep=1pt,draw=black, fill=black]
\begin{tikzpicture}
\node[point] (h0) at (0.5773502691896257, 1.0) {};
\node[point] (h2) at (0.28867513459481287, 0.5) {};
\node[point] (h4) at (1.4433756729740643, 0.5) {};
\node[point] (h5) at (1.1547005383792515, 1.0) {};
\node[point] (h6) at (0.5773502691896257, 0.0) {};
\node[point] (h7) at (1.1547005383792515, 0.0) {};
\draw (0.28867513459481287, 0.5) -- (0.5773502691896257, 0.0);
\draw (0.28867513459481287, 0.5) -- (0.5773502691896257, 1.0);
\draw (0.5773502691896257, 0.0) -- (1.1547005383792515, 0.0);
\draw (0.5773502691896257, 1.0) -- (1.1547005383792515, 1.0);
\draw (1.1547005383792515, 0.0) -- (1.4433756729740643, 0.5);
\draw (1.1547005383792515, 1.0) -- (1.4433756729740643, 0.5);
\end{tikzpicture}\hspace{15mm} \tikzstyle{point}=[circle,minimum size=1pt, inner sep=1pt,draw=black, fill=black]
\begin{tikzpicture}
\node[point] (h0) at (2.309401076758503, 2.0) {};
\node[point] (h1) at (2.0207259421636903, 2.5) {};
\node[point] (h2) at (0.5773502691896257, 1.0) {};
\node[point] (h3) at (1.4433756729740643, 1.5) {};
\node[point] (h4) at (0.28867513459481287, 0.5) {};
\node[point] (h5) at (3.1754264805429417, 0.5) {};
\node[point] (h6) at (1.1547005383792515, 0.0) {};
\node[point] (h7) at (2.8867513459481287, 0.0) {};
\node[point] (h8) at (2.309401076758503, 1.0) {};
\node[point] (h9) at (2.0207259421636903, 1.5) {};
\node[point] (h10) at (1.4433756729740643, -0.5) {};
\node[point] (h11) at (1.4433756729740643, 0.5) {};
\node[point] (h12) at (0.5773502691896257, 0.0) {};
\node[point] (h13) at (2.8867513459481287, 2.0) {};
\node[point] (h14) at (1.1547005383792515, 2.0) {};
\node[point] (h15) at (2.309401076758503, 0.0) {};
\node[point] (h16) at (2.0207259421636903, 0.5) {};
\node[point] (h17) at (2.0207259421636903, -0.5) {};
\node[point] (h18) at (0.5773502691896257, 2.0) {};
\node[point] (h19) at (1.4433756729740643, 2.5) {};
\node[point] (h20) at (3.1754264805429417, 1.5) {};
\node[point] (h21) at (0.28867513459481287, 1.5) {};
\node[point] (h22) at (1.1547005383792515, 1.0) {};
\node[point] (h23) at (2.8867513459481287, 1.0) {};
\draw (0.28867513459481287, 0.5) -- (0.5773502691896257, 0.0);
\draw (0.28867513459481287, 0.5) -- (0.5773502691896257, 1.0);
\draw (0.5773502691896257, 0.0) -- (1.1547005383792515, 0.0);
\draw (0.28867513459481287, 1.5) -- (0.5773502691896257, 1.0);
\draw (0.28867513459481287, 1.5) -- (0.5773502691896257, 2.0);
\draw (0.5773502691896257, 1.0) -- (1.1547005383792515, 1.0);
\draw (0.5773502691896257, 2.0) -- (1.1547005383792515, 2.0);
\draw (1.1547005383792515, 0.0) -- (1.4433756729740643, -0.5);
\draw (1.1547005383792515, 0.0) -- (1.4433756729740643, 0.5);
\draw (1.4433756729740643, -0.5) -- (2.0207259421636903, -0.5);
\draw (1.1547005383792515, 1.0) -- (1.4433756729740643, 0.5);
\draw (1.1547005383792515, 1.0) -- (1.4433756729740643, 1.5);
\draw (1.4433756729740643, 0.5) -- (2.0207259421636903, 0.5);
\draw (1.1547005383792515, 2.0) -- (1.4433756729740643, 1.5);
\draw (1.1547005383792515, 2.0) -- (1.4433756729740643, 2.5);
\draw (1.4433756729740643, 1.5) -- (2.0207259421636903, 1.5);
\draw (1.4433756729740643, 2.5) -- (2.0207259421636903, 2.5);
\draw (2.0207259421636903, -0.5) -- (2.309401076758503, 0.0);
\draw (2.0207259421636903, 0.5) -- (2.309401076758503, 0.0);
\draw (2.0207259421636903, 0.5) -- (2.309401076758503, 1.0);
\draw (2.309401076758503, 0.0) -- (2.8867513459481287, 0.0);
\draw (2.0207259421636903, 1.5) -- (2.309401076758503, 1.0);
\draw (2.0207259421636903, 1.5) -- (2.309401076758503, 2.0);
\draw (2.309401076758503, 1.0) -- (2.8867513459481287, 1.0);
\draw (2.0207259421636903, 2.5) -- (2.309401076758503, 2.0);
\draw (2.309401076758503, 2.0) -- (2.8867513459481287, 2.0);
\draw (2.8867513459481287, 0.0) -- (3.1754264805429417, 0.5);
\draw (2.8867513459481287, 1.0) -- (3.1754264805429417, 0.5);
\draw (2.8867513459481287, 1.0) -- (3.1754264805429417, 1.5);
\draw (2.8867513459481287, 2.0) -- (3.1754264805429417, 1.5);
\end{tikzpicture}\hspace{15mm} \tikzstyle{point}=[circle,minimum size=1pt, inner sep=1pt,draw=black, fill=black]
\begin{tikzpicture}
\node[point] (h0) at (2.309401076758503, 2.0) {};
\node[point] (h1) at (2.0207259421636903, 2.5) {};
\node[point] (h2) at (4.618802153517006, 1.0) {};
\node[point] (h3) at (0.5773502691896257, 1.0) {};
\node[point] (h4) at (1.4433756729740643, 1.5) {};
\node[point] (h5) at (4.04145188432738, 0.0) {};
\node[point] (h6) at (0.28867513459481287, 0.5) {};
\node[point] (h7) at (3.1754264805429417, 3.5) {};
\node[point] (h8) at (4.907477288111819, 0.5) {};
\node[point] (h9) at (4.04145188432738, 3.0) {};
\node[point] (h10) at (3.1754264805429417, -0.5) {};
\node[point] (h11) at (3.1754264805429417, 0.5) {};
\node[point] (h12) at (1.1547005383792515, 0.0) {};
\node[point] (h13) at (2.8867513459481287, 0.0) {};
\node[point] (h14) at (1.1547005383792515, 3.0) {};
\node[point] (h15) at (2.8867513459481287, 3.0) {};
\node[point] (h16) at (3.752776749732567, 2.5) {};
\node[point] (h17) at (2.309401076758503, 1.0) {};
\node[point] (h18) at (2.0207259421636903, 1.5) {};
\node[point] (h19) at (2.309401076758503, -1.0) {};
\node[point] (h20) at (1.4433756729740643, -0.5) {};
\node[point] (h21) at (1.4433756729740643, 0.5) {};
\node[point] (h22) at (2.309401076758503, 4.0) {};
\node[point] (h23) at (4.618802153517006, 0.0) {};
\node[point] (h24) at (0.5773502691896257, 3.0) {};
\node[point] (h25) at (1.4433756729740643, 3.5) {};
\node[point] (h26) at (4.04145188432738, 2.0) {};
\node[point] (h27) at (0.5773502691896257, 0.0) {};
\node[point] (h28) at (4.618802153517006, 3.0) {};
\node[point] (h29) at (3.1754264805429417, 2.5) {};
\node[point] (h30) at (0.28867513459481287, 2.5) {};
\node[point] (h31) at (4.907477288111819, 2.5) {};
\node[point] (h32) at (2.8867513459481287, 2.0) {};
\node[point] (h33) at (3.752776749732567, 1.5) {};
\node[point] (h34) at (1.1547005383792515, 2.0) {};
\node[point] (h35) at (2.309401076758503, 3.0) {};
\node[point] (h36) at (2.0207259421636903, 3.5) {};
\node[point] (h37) at (2.309401076758503, 0.0) {};
\node[point] (h38) at (2.0207259421636903, 0.5) {};
\node[point] (h39) at (2.0207259421636903, -0.5) {};
\node[point] (h40) at (4.618802153517006, 2.0) {};
\node[point] (h41) at (3.752776749732567, 0.5) {};
\node[point] (h42) at (0.5773502691896257, 2.0) {};
\node[point] (h43) at (1.4433756729740643, 2.5) {};
\node[point] (h44) at (3.1754264805429417, 1.5) {};
\node[point] (h45) at (0.28867513459481287, 1.5) {};
\node[point] (h46) at (4.04145188432738, 1.0) {};
\node[point] (h47) at (4.907477288111819, 1.5) {};
\node[point] (h48) at (2.8867513459481287, -1.0) {};
\node[point] (h49) at (1.1547005383792515, 1.0) {};
\node[point] (h50) at (2.8867513459481287, 4.0) {};
\node[point] (h51) at (2.8867513459481287, 1.0) {};
\node[point] (h52) at (3.752776749732567, -0.5) {};
\node[point] (h53) at (3.752776749732567, 3.5) {};
\draw (0.28867513459481287, 0.5) -- (0.5773502691896257, 0.0);
\draw (0.28867513459481287, 0.5) -- (0.5773502691896257, 1.0);
\draw (0.5773502691896257, 0.0) -- (1.1547005383792515, 0.0);
\draw (0.28867513459481287, 1.5) -- (0.5773502691896257, 1.0);
\draw (0.28867513459481287, 1.5) -- (0.5773502691896257, 2.0);
\draw (0.5773502691896257, 1.0) -- (1.1547005383792515, 1.0);
\draw (0.28867513459481287, 2.5) -- (0.5773502691896257, 2.0);
\draw (0.28867513459481287, 2.5) -- (0.5773502691896257, 3.0);
\draw (0.5773502691896257, 2.0) -- (1.1547005383792515, 2.0);
\draw (0.5773502691896257, 3.0) -- (1.1547005383792515, 3.0);
\draw (1.1547005383792515, 0.0) -- (1.4433756729740643, -0.5);
\draw (1.1547005383792515, 0.0) -- (1.4433756729740643, 0.5);
\draw (1.4433756729740643, -0.5) -- (2.0207259421636903, -0.5);
\draw (1.1547005383792515, 1.0) -- (1.4433756729740643, 0.5);
\draw (1.1547005383792515, 1.0) -- (1.4433756729740643, 1.5);
\draw (1.4433756729740643, 0.5) -- (2.0207259421636903, 0.5);
\draw (1.1547005383792515, 2.0) -- (1.4433756729740643, 1.5);
\draw (1.1547005383792515, 2.0) -- (1.4433756729740643, 2.5);
\draw (1.4433756729740643, 1.5) -- (2.0207259421636903, 1.5);
\draw (1.1547005383792515, 3.0) -- (1.4433756729740643, 2.5);
\draw (1.1547005383792515, 3.0) -- (1.4433756729740643, 3.5);
\draw (1.4433756729740643, 2.5) -- (2.0207259421636903, 2.5);
\draw (1.4433756729740643, 3.5) -- (2.0207259421636903, 3.5);
\draw (2.0207259421636903, -0.5) -- (2.309401076758503, -1.0);
\draw (2.0207259421636903, -0.5) -- (2.309401076758503, 0.0);
\draw (2.309401076758503, -1.0) -- (2.8867513459481287, -1.0);
\draw (2.0207259421636903, 0.5) -- (2.309401076758503, 0.0);
\draw (2.0207259421636903, 0.5) -- (2.309401076758503, 1.0);
\draw (2.309401076758503, 0.0) -- (2.8867513459481287, 0.0);
\draw (2.0207259421636903, 1.5) -- (2.309401076758503, 1.0);
\draw (2.0207259421636903, 1.5) -- (2.309401076758503, 2.0);
\draw (2.309401076758503, 1.0) -- (2.8867513459481287, 1.0);
\draw (2.0207259421636903, 2.5) -- (2.309401076758503, 2.0);
\draw (2.0207259421636903, 2.5) -- (2.309401076758503, 3.0);
\draw (2.309401076758503, 2.0) -- (2.8867513459481287, 2.0);
\draw (2.0207259421636903, 3.5) -- (2.309401076758503, 3.0);
\draw (2.0207259421636903, 3.5) -- (2.309401076758503, 4.0);
\draw (2.309401076758503, 3.0) -- (2.8867513459481287, 3.0);
\draw (2.309401076758503, 4.0) -- (2.8867513459481287, 4.0);
\draw (2.8867513459481287, -1.0) -- (3.1754264805429417, -0.5);
\draw (2.8867513459481287, 0.0) -- (3.1754264805429417, -0.5);
\draw (2.8867513459481287, 0.0) -- (3.1754264805429417, 0.5);
\draw (3.1754264805429417, -0.5) -- (3.752776749732567, -0.5);
\draw (2.8867513459481287, 1.0) -- (3.1754264805429417, 0.5);
\draw (2.8867513459481287, 1.0) -- (3.1754264805429417, 1.5);
\draw (3.1754264805429417, 0.5) -- (3.752776749732567, 0.5);
\draw (2.8867513459481287, 2.0) -- (3.1754264805429417, 1.5);
\draw (2.8867513459481287, 2.0) -- (3.1754264805429417, 2.5);
\draw (3.1754264805429417, 1.5) -- (3.752776749732567, 1.5);
\draw (2.8867513459481287, 3.0) -- (3.1754264805429417, 2.5);
\draw (2.8867513459481287, 3.0) -- (3.1754264805429417, 3.5);
\draw (3.1754264805429417, 2.5) -- (3.752776749732567, 2.5);
\draw (2.8867513459481287, 4.0) -- (3.1754264805429417, 3.5);
\draw (3.1754264805429417, 3.5) -- (3.752776749732567, 3.5);
\draw (3.752776749732567, -0.5) -- (4.04145188432738, 0.0);
\draw (3.752776749732567, 0.5) -- (4.04145188432738, 0.0);
\draw (3.752776749732567, 0.5) -- (4.04145188432738, 1.0);
\draw (4.04145188432738, 0.0) -- (4.618802153517006, 0.0);
\draw (3.752776749732567, 1.5) -- (4.04145188432738, 1.0);
\draw (3.752776749732567, 1.5) -- (4.04145188432738, 2.0);
\draw (4.04145188432738, 1.0) -- (4.618802153517006, 1.0);
\draw (3.752776749732567, 2.5) -- (4.04145188432738, 2.0);
\draw (3.752776749732567, 2.5) -- (4.04145188432738, 3.0);
\draw (4.04145188432738, 2.0) -- (4.618802153517006, 2.0);
\draw (3.752776749732567, 3.5) -- (4.04145188432738, 3.0);
\draw (4.04145188432738, 3.0) -- (4.618802153517006, 3.0);
\draw (4.618802153517006, 0.0) -- (4.907477288111819, 0.5);
\draw (4.618802153517006, 1.0) -- (4.907477288111819, 0.5);
\draw (4.618802153517006, 1.0) -- (4.907477288111819, 1.5);
\draw (4.618802153517006, 2.0) -- (4.907477288111819, 1.5);
\draw (4.618802153517006, 2.0) -- (4.907477288111819, 2.5);
\draw (4.618802153517006, 3.0) -- (4.907477288111819, 2.5);
\end{tikzpicture}
\caption[LoF entry]{Visualization of the honeycomb networks $HC(1)$, $HC(2)$, and $HC(3)$ from left to right.}
\label{fig:honeycomb_examples}
\end{figure}
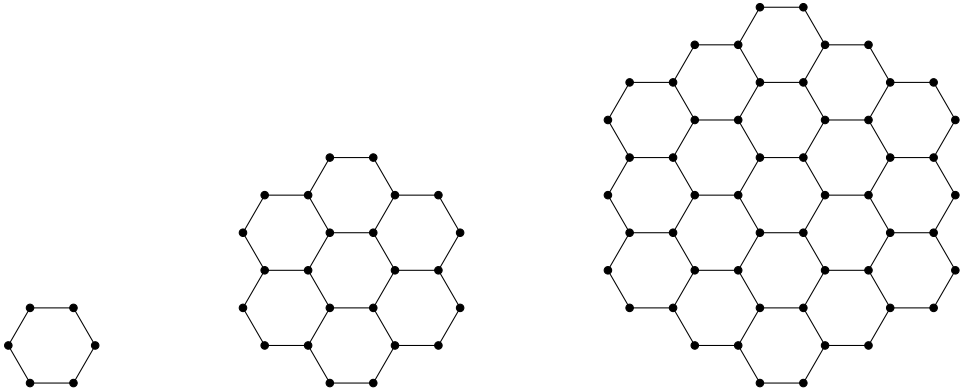

Naturally, two-dimensional grids correspond to a square tiling of $\mathbb{R}^2$. Similarly, following the terminology and notation of ~\cite{manuel2008minimum}, the honeycomb networks $HC(n)$ correspond to a partial hexagonal tiling of $\mathbb{R}^2$.
More precisely, $HC(1)$ is a single regular hexagon and $HC(n)$ consists of $n-1$ layers of hexagons around a central one (see~\cref{fig:honeycomb_examples}, left). The number of vertices in $HC(n)$ is $6n^2$.

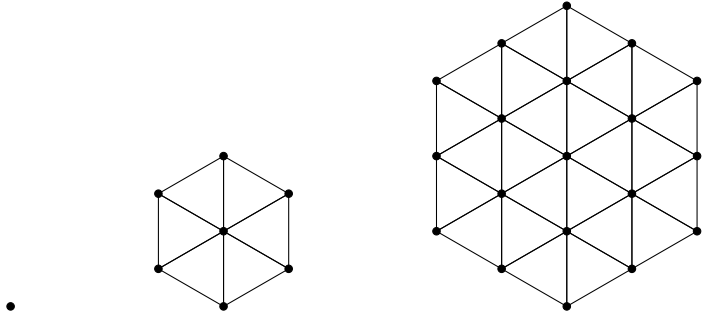
\begin{figure}[h]
\centering
\tikzstyle{point}=[circle,minimum size=1pt, inner sep=1pt,draw=black, fill=black]
\begin{tikzpicture}
\node[point] (1) at (0.0, 0) {};
\end{tikzpicture}\hspace{15mm} \tikzstyle{point}=[circle,minimum size=1pt, inner sep=1pt,draw=black, fill=black]
\begin{tikzpicture}
\node[point] (1) at (0.0, 0) {};
\node[point] (2) at (0.0, 1) {};
\node[point] (3) at (0.8660254037844386, -0.5) {};
\node[point] (4) at (0.8660254037844386, 0.5) {};
\node[point] (5) at (0.8660254037844386, 1.5) {};
\node[point] (7) at (1.7320508075688772, 0.0) {};
\node[point] (8) at (1.7320508075688772, 1.0) {};
\draw (0.0, 0) -- (0.0, 1) -- (0.8660254037844386, 0.5) -- (0.0, 0);
\draw (0.0, 0) -- (0.8660254037844386, -0.5) -- (0.8660254037844386, 0.5) -- (0.0, 0);
\draw (0.0, 1) -- (0.8660254037844386, 0.5) -- (0.8660254037844386, 1.5) -- (0.0, 1);
\draw (0.8660254037844386, -0.5) -- (0.8660254037844386, 0.5) -- (1.7320508075688772, 0.0) -- (0.8660254037844386, -0.5);
\draw (0.8660254037844386, 0.5) -- (0.8660254037844386, 1.5) -- (1.7320508075688772, 1.0) -- (0.8660254037844386, 0.5);
\draw (0.8660254037844386, 0.5) -- (1.7320508075688772, 0.0) -- (1.7320508075688772, 1.0) -- (0.8660254037844386, 0.5);
\end{tikzpicture}\hspace{15mm} \tikzstyle{point}=[circle,minimum size=1pt, inner sep=1pt,draw=black, fill=black]
\begin{tikzpicture}
\node[point] (1) at (0.0, 0) {};
\node[point] (2) at (0.0, 1) {};
\node[point] (3) at (0.0, 2) {};
\node[point] (4) at (0.8660254037844386, -0.5) {};
\node[point] (5) at (0.8660254037844386, 0.5) {};
\node[point] (6) at (0.8660254037844386, 1.5) {};
\node[point] (7) at (0.8660254037844386, 2.5) {};
\node[point] (8) at (1.7320508075688772, -1.0) {};
\node[point] (9) at (1.7320508075688772, 0.0) {};
\node[point] (10) at (1.7320508075688772, 1.0) {};
\node[point] (11) at (1.7320508075688772, 2.0) {};
\node[point] (12) at (1.7320508075688772, 3.0) {};
\node[point] (13) at (2.598076211353316, -0.5) {};
\node[point] (14) at (2.598076211353316, 0.5) {};
\node[point] (15) at (2.598076211353316, 1.5) {};
\node[point] (16) at (2.598076211353316, 2.5) {};
\node[point] (17) at (3.4641016151377544, 0.0) {};
\node[point] (18) at (3.4641016151377544, 1.0) {};
\node[point] (19) at (3.4641016151377544, 2.0) {};
\draw (0.0, 0) -- (0.0, 1) -- (0.8660254037844386, 0.5) -- (0.0, 0);
\draw (0.0, 0) -- (0.8660254037844386, -0.5) -- (0.8660254037844386, 0.5) -- (0.0, 0);
\draw (0.0, 1) -- (0.0, 2) -- (0.8660254037844386, 1.5) -- (0.0, 1);
\draw (0.0, 1) -- (0.8660254037844386, 0.5) -- (0.8660254037844386, 1.5) -- (0.0, 1);
\draw (0.0, 2) -- (0.8660254037844386, 1.5) -- (0.8660254037844386, 2.5) -- (0.0, 2);
\draw (0.8660254037844386, -0.5) -- (0.8660254037844386, 0.5) -- (1.7320508075688772, 0.0) -- (0.8660254037844386, -0.5);
\draw (0.8660254037844386, -0.5) -- (1.7320508075688772, -1.0) -- (1.7320508075688772, 0.0) -- (0.8660254037844386, -0.5);
\draw (0.8660254037844386, 0.5) -- (0.8660254037844386, 1.5) -- (1.7320508075688772, 1.0) -- (0.8660254037844386, 0.5);
\draw (0.8660254037844386, 0.5) -- (1.7320508075688772, 0.0) -- (1.7320508075688772, 1.0) -- (0.8660254037844386, 0.5);
\draw (0.8660254037844386, 1.5) -- (0.8660254037844386, 2.5) -- (1.7320508075688772, 2.0) -- (0.8660254037844386, 1.5);
\draw (0.8660254037844386, 1.5) -- (1.7320508075688772, 1.0) -- (1.7320508075688772, 2.0) -- (0.8660254037844386, 1.5);
\draw (0.8660254037844386, 2.5) -- (1.7320508075688772, 2.0) -- (1.7320508075688772, 3.0) -- (0.8660254037844386, 2.5);
\draw (1.7320508075688772, -1.0) -- (1.7320508075688772, 0.0) -- (2.598076211353316, -0.5) -- (1.7320508075688772, -1.0);
\draw (1.7320508075688772, 0.0) -- (1.7320508075688772, 1.0) -- (2.598076211353316, 0.5) -- (1.7320508075688772, 0.0);
\draw (1.7320508075688772, 0.0) -- (2.598076211353316, -0.5) -- (2.598076211353316, 0.5) -- (1.7320508075688772, 0.0);
\draw (1.7320508075688772, 1.0) -- (1.7320508075688772, 2.0) -- (2.598076211353316, 1.5) -- (1.7320508075688772, 1.0);
\draw (1.7320508075688772, 1.0) -- (2.598076211353316, 0.5) -- (2.598076211353316, 1.5) -- (1.7320508075688772, 1.0);
\draw (1.7320508075688772, 2.0) -- (1.7320508075688772, 3.0) -- (2.598076211353316, 2.5) -- (1.7320508075688772, 2.0);
\draw (1.7320508075688772, 2.0) -- (2.598076211353316, 1.5) -- (2.598076211353316, 2.5) -- (1.7320508075688772, 2.0);
\draw (2.598076211353316, -0.5) -- (2.598076211353316, 0.5) -- (3.4641016151377544, 0.0) -- (2.598076211353316, -0.5);
\draw (2.598076211353316, 0.5) -- (2.598076211353316, 1.5) -- (3.4641016151377544, 1.0) -- (2.598076211353316, 0.5);
\draw (2.598076211353316, 0.5) -- (3.4641016151377544, 0.0) -- (3.4641016151377544, 1.0) -- (2.598076211353316, 0.5);
\draw (2.598076211353316, 1.5) -- (2.598076211353316, 2.5) -- (3.4641016151377544, 2.0) -- (2.598076211353316, 1.5);
\draw (2.598076211353316, 1.5) -- (3.4641016151377544, 1.0) -- (3.4641016151377544, 2.0) -- (2.598076211353316, 1.5);
\end{tikzpicture}
\caption[LoF entry]{Visualization of the hexagon networks $HX(1)$, $HX(2)$, and $HX(3)$ from left to right.}
\label{fig:hex_examples}
\end{figure}

The hexagon networks $HX(n)$ correspond also to a recursive and partial tiling of $\mathbb{R}^2$ but with equilateral triangles. $HX(1)$ is, by definition, a single vertex. $HX(n)$ is a regular hexagon with sides of length $s$ tiled with equilateral triangles with sides of length $\frac{s}{n-1}$. In particular, $HX(n)$ has $3n^2-3n+1$ vertices from $6(n-1)^2$ triangles (see Figure~\ref{fig:hex_examples}).

\begin{figure}[h]
\centering
\input{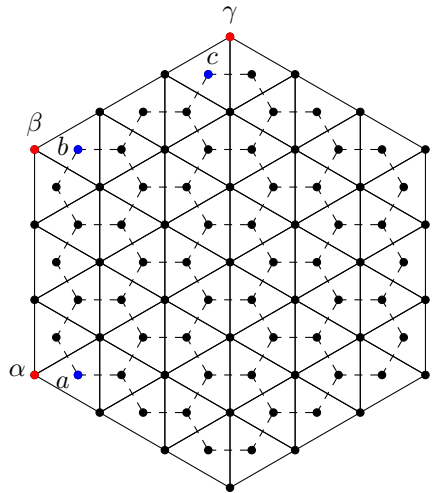}
\caption[LoF entry]{Visualization of the honeycomb network $HC(3)$ (dashed) as the bounded dual of the hexagonal network $HX(4)$ with minimum resolving sets $\{a,b,c\}$ and $\{\alpha,\beta,\gamma\}$, respectively, following the construction in~\cite{manuel2008minimum}.}
\label{fig:honeycomb_hex_res_set_ex}
\end{figure}

By identifying a useful coordinate system for dealing with distances in hexagonal graphs, one can show that any three vertices of degree three forming adjacent corners of the underlying hexagon in $HX(n)$ serve as a resolving set (see~\cref{fig:honeycomb_hex_res_set_ex}). Since $HX(n)$ violates properties that any graph with metric dimension two must have~\cite{khuller1996landmarks}, we conclude that $\beta(HX(n)) = 3$~\cite{manuel2008minimum}.

To determine the metric dimension of $HC(n)$, it is useful to note that $HC(n)$ is the so-called bounded dual of $HX(n)$. In particular, $HC(n)$ may be constructed from $HX(n+1)$ as follows. For each face of $HX(n+1)$, except the unbounded face, include a single vertex in $HC(n)$. Two vertices in $HC(n)$ are adjacent when the corresponding faces in $HX(n+1)$ share an edge. By taking advantage of this relationship, it can be shown that $\beta(HC(n)) = 3$ as well~\cite{manuel2008minimum} (see~\cref{fig:honeycomb_hex_res_set_ex}). 

\subsection{Hamming Graphs}
\label{sec:hamming}
There are many ways to measure distance between pairs of strings. One of the simplest and most common is the Hamming distance~\cite{hamming1950error}. Comparing two strings of the same length, the Hamming distance counts the number of positions in which the strings disagree. This distance can be used to define a graph.

\begin{definition}{(Hamming Graph.)}
Let $V = A^k$ be the set of all strings of length $k$ from the alphabet $A$ of size $a=|A|$, and let $d(u,v)$ be the Hamming distance between $u,v \in V$. The Hamming graph $H_{k,a}$ contains a vertex associated with each $v \in V$ and the edge $\{u,v\} \in E$ only when $d(u,v)=1$. 
\end{definition}

Despite the highly symmetric nature of Hamming graphs, an efficient algorithm to compute their metric dimension is not known. 
However, some upper bounds on $\beta(H_{k, a})$ are known.

For the special case when $a = 2$, $Q_k := H_{k,2}$ is called a hypercube.
The hypercube turns out to be an important structure when studying a particular coin weighing problem: given a set of $k$ coins of two different weights, how many weighings are required to determine which coins are lighter~\cite{soderberg1963combinatory}? An asymptotic result tells us that $\lim_{k \rightarrow \infty} \beta(Q_k) \frac{\log_2(k)}{k} = 2$, suggesting a method using $2^n-1$ weighings to find the answer when $k = n2^{n-1}$~\cite{erdos1963two,lindstrom1964combinatory}.
This asymptotic behavior generalizes to arbitrary Hamming graphs as $\lim_{k \rightarrow \infty} \beta(H_{k,a}) \frac{\log_a(k)}{k} = 2$~\cite{jiang2019metric}.

Analysis of the Mastermind game, a game closely related to the coin weighing problem above, depends on the structure of Hamming graphs as well. Mastermind is played between two players, $A$ and $B$. Player $A$ starts by choosing $s = s_1\dots s_k$, a sequence unknown to player $B$ of $k$ symbols from an alphabet of size $c$. Player $B$ attempts to guess $s$ as quickly possible. After each guess $q$, player $A$ provides two values: $a(q, s)$, the number of positions where $q$ and $s$ agree, and $b(q, s)$, the total number of correct symbols at incorrect positions. Donald Knuth created an algorithm solving the commercial version of the game with $k=4$ and $c=6$ with at most five questions~\cite{knuth1976computer}. For arbitrary values of $k$ and $c$ a number of results and bounds exist~\cite{chvatal1983mastermind,goddard2004mastermind} though there is no known optimal solution. 

A static variant of the game in which player $B$ must make all guesses at once with no feedback and player $A$ only provides $a(q,s)$ has been shown to be NP-complete~\cite{goodrich2009algorithmic}. Indeed, notice that $a(q, s) = k - d(q,s)$, where $d(q,s)$ is the Hamming distance between the two sequences. In particular, to guess $s$ with the least number of questions, player $B$ should make guesses that reveal the Hamming distance from $s$ to a resolving set of $H_{k,c}$ of size $\beta(H_{k,c})$. A deeper analysis of this situation yields an upper bound on $\beta(H_{k,c})$ when $k$ is small in comparison to $c$. In particular, if $\epsilon < 1$ and $c < k^{1-\epsilon}$ then $\beta(H_{k,c}) \leq (2+\epsilon)k\frac{1+2\log_2(c)}{\log_2(k)-\log_2(c)}$~\cite{chvatal1983mastermind}.

Bounds related to the metric dimension of Cartesian product of graphs have proven useful in studying Hamming graphs.
In particular, $\beta(G) \leq \beta(G \square K_2) \leq \beta(G)+1$~\cite{chartrand2000resolvability} and, since $H_{k,a} = K_a^{\square k}$ is the Cartesian product of $k$ copies of the complete graph $K_a$, it follows that $\beta(Q_k) \leq \beta(Q_{k+1}) \leq \beta(Q_k) + 1$.
This approach yields the identity $\beta(H_{2,a}) = \lfloor \frac{2}{3} (2a-1) \rfloor$~\cite{caceres2007metric}. 
\cref{tab:hypercube_betas} shows exact values of $\beta(Q_k)$ verified via brute force for $1 \leq k \leq 10$, and upper bounds based on a variable neighborhood search for $11 \leq k \leq 17$~\cite{mladenovic2012variable}.

\begin{table}[h!]
\centering 
\small
\begin{tabular}{ccccccccccc|ccccccc}
\toprule
$k$ & 1 & 2 & 3 & 4 & 5 & 6 & 7 & 8 & 9 & 10 & 11 & 12 & 13 & 14 & 15 & 16 & 17\\
$\beta(Q_k)$ & 1 & 2 & 3 & 4 & 4 & 5 & 6 & 6 & 7 & 7 & 8 & 8 & 8 & 9 & 9 & 10 & 10 \\ 
\bottomrule\\
\end{tabular}
\caption{Exact values of $\beta(Q_k)$ for $1 \leq k \leq 10$, and upper bounds for $11 \leq k \leq 17$~\cite{mladenovic2012variable}.}
\label{tab:hypercube_betas}
\end{table}

The bounds for the metric dimension of hypercubes have been generalized for arbitrary Hamming graphs as $\beta(H_{k,a}) \leq \beta(H_{k+1,a}) \leq \beta(H_{k,a}) + \lfloor \frac{a}{2} \rfloor$~\cite{TilLla18}. The proof is constructive allowing for the generation of $R_{k+1}$, a resolving set for $H_{k+1,a}$, from $R_k$, any resolving set for $H_{k,a}$. 

When dealing with very large Hamming graphs, verifying that a given subset of vertices is resolving becomes intractable via the simple brute force approach. However, the recursive structure and highly symmetric nature of Hamming graphs allows for far more efficient resolvability checks in practice. Indeed, by describing resolvability on Hamming graphs as a linear system, integer programming techniques can be used to verify quickly that a given set of vertices is resolving with high probability. A somewhat slower but deterministic solution can be implemented using Gr\"obner bases. These techniques have been used to discover a resolving set of size 77 for the Hamming graph $H_{8,20}$, showing that $\beta(H_{8,20}) \leq 77$~\cite{laird2019resolvability}. With 25.6 billion vertices, traditional methods for finding small resolving sets are not computationally feasible in this setting.

\begin{table}
\centering 
\scriptsize
\begin{tabular}{llll}
\toprule
\textbf{Graph Type} & \textbf{Symbol (constraint)} & \textbf{Metric Dimension} & \textbf{Reference} \\
\cmidrule(lr){1-1}
\cmidrule(lr){2-2}
\cmidrule(lr){3-3}
\cmidrule(lr){4-4}
Antiprisms & $A_n\,(n \geq 3$) & 3 &  \cite{javaid2008families} \\ 
Complete Graphs & $K_n$ & $n-1$ &  \cite{chartrand2000resolvability} \\
Cycles & $C_n$ & 2 &  \cite{chartrand2000resolvabilityCycles} \\
De Bruijn Graphs & $B_{d, n}$ & $d^{n-1}(d-1)$ & \cite{feng2013metric,rajan2014metric} \\ 
Fans & $F_n\,(n \not \in \{1,2,3,6\})$ & $\lfloor \frac{2n+2}{5} \rfloor$ &  \cite{hernando2005metric} \\ 
Grids in $d$ dimensions & $G_{m,n,\dots}$ & $d$ &  \cite{khuller1996landmarks,melter1984metric} \\ 
Hexagonal Graphs & $HX(n)$ & 3 &  \cite{manuel2008minimum} \\ 
Honeycomb Graphs & $HC(n)$ & 3 &  \cite{manuel2008minimum} \\ 
Jahangir Graphs & $J_{2n}\,(n \geq 4)$ & $\lfloor \frac{2n}{3} \rfloor$ &  \cite{tomescu2007metric} \\ 
Kautz Graphs & $K_{d, n}$ & $(d^{n-1} + d^{n-2})(d-1)$ & \cite{feng2013metric,rajan2014metric} \\ 
Paths & $P_n$ & 1 &  \cite{chartrand2000resolvability}\\ 
Peterson Graphs & $P_{n, 2}\,(n \geq 5)$ & 3 &  \cite{javaid2008families} \\ 
Prisms & $D_n$ & 2 for odd $n$, 3 for even $n$ &  \cite{javaid2008families} \\ 
Trees & $T_n$ & $\ell(T_n) - \text{ex}(T_n)$  &  \cite{chartrand2000resolvability,harary1976metric,slater1975leaves} \\ 
Wheels & $W_n\,(n>6)$ & $\beta(W_n) = \beta(F_n)$ &  \S~\ref{sec:fans-wheels}, \cite{shanmukha2002metric} \\
\bottomrule\\
\end{tabular}
\caption{Exact values for the metric dimension of several different families of graphs.}
\label{tab:beta_list_exact}
\end{table}

\begin{table}
\centering 
\tiny
\begin{tabular}{llll}
\toprule
\textbf{Graph Type} & \textbf{Symbol (constraint)} & \textbf{Metric Dimension} & \textbf{Reference} \\
\cmidrule(lr){1-1}
\cmidrule(lr){2-2}
\cmidrule(lr){3-3}
\cmidrule(lr){4-4}
Bilinear Forms Graphs & $H_q(n,d)\,(n \geq d \geq 2)$ & $\beta(H_q(n,d)) \leq q^{n+d-1+ \lfloor \frac{d+1}{n} \rfloor}$ &  \cite{feng2012metric} \\
Cartesian Product with $K_2$ & $G \times K_2$ & $\beta(G) \leq \beta(G \times K_2) \leq \beta(G)+1$ &  \cite{chartrand2000resolvability} \\
Doubled Odd Graphs & $O(2e+1,e,e+1)\,(e \geq 2)$ & $\beta(O(2e+1,e,e+1)) \leq 2e+1$ &  \cite{guo2013metric} \\
Hamming Graphs & $H_{k, a}$ & $\beta(H_{k,a}) \leq \beta(H_{k+1, a}) \leq \beta(H_{k, a})\!+\!\lfloor a/2\rfloor$ &  \cite{TilLla18} \\
Johnson Graphs & $J(n,e)$ & $\beta(J(n,e)) \leq (e+1)\lceil n/(e+1) \rceil$ & \cite{bailey2013resolving,guo2013metric} \\
Unicyclic Graphs & $T+e$ & $\beta(T) - 2 \leq \beta(T + e) \leq \beta(T)+1$ &  \cite{chartrand2000resolvability} \\ 
\bottomrule\\
\end{tabular}
\caption{Bounds on the metric dimension of several different families of graphs.}
\label{tab:beta_list_bounds}
\end{table}

\begin{table}
\centering 
\scriptsize
\begin{tabular}{ll}
\toprule
\textbf{Graph Type} & \textbf{Metric Dimension} \\
\cmidrule(lr){1-1}
\cmidrule(lr){2-2}
Benes Networks & Polynomial time solvable, see reference \cite{manuel2008efficient} \\ 
Butterfly Networks & Polynomial time solvable, see reference \cite{manuel2008efficient} \\ 
Cactus Block Graphs & Linear time solvable, see reference \cite{hoffmann2016linear} \\ 
Chain Graphs & Linear time solvable, see reference \cite{fernau2015computing} \\ 
Cographs & Linear time solvable, see reference \cite{epstein2015weighted} \\ 
Outerplanar Graphs & Polynomial time solvable, see reference \cite{diaz2012complexity} \\ 
\bottomrule\\
\end{tabular}
\caption{Families of graphs for which metric dimension can be determined efficiently.}
\label{tab:beta_list_algo}
\end{table}

\begin{table}
\centering 
\scriptsize
\begin{tabular}{ll}
\toprule
\textbf{Graph Type} & \textbf{Metric Dimension} \\
\cmidrule(lr){1-1}
\cmidrule(lr){2-2}
Amalgamation of Cycles & See page 25 of reference \cite{iswadi2010metric} \\ 
Cayley Digraphs & See pages 34-37 of reference \cite{fehr2006metric} \\ 
Circulent Networks & See reference \cite{rajan2010minimum} \\  
Complete $k$-partite Graphs & See reference \cite{saputro2009metric} \\ 
Generalized Wheel Graphs & See reference \cite{sooryanarayana2019metric} \\ 
Graphs with Pendant Edges & See pages 4-7 of reference \cite{iswadi2008metric} \\ 
Grassmann Graphs & See page 98 of reference \cite{meagher2012metric} \\ 
Harary Graphs $H_{4,n}$ & See page 9 of reference \cite{javaid2008families} \\ 
Kneser Graphs & See page 750 of reference \cite{bailey2013resolving} \\ 
Line Graphs & See pages 803-804 of reference \cite{feng2013metric} \\ 
Regular Bipartite Graphs & See pages 16-17 of reference \cite{bavca2011metric} \\ 
Torus Network & See pages 268 and 271 of reference \cite{manuel2006landmarks} \\ 
Twisted Grassmann Graphs & See page 4 of reference \cite{guo2013metric} \\
\bottomrule\\
 \end{tabular}
\caption{Additional families of graphs for which metric dimension can be determined efficiently. See \cite{hernando2005metric} for additional bounds.}
 \label{tab:beta_list_other}
\end{table}

\section{Random Graph Models}

Real world networks rarely fully conform to the requirements for structurally deterministic graph families.
Random graph models, which define distributions over graph structures, often in terms of some generative process, are therefore useful in describing real networks.
Understanding the behavior of metric dimension as a random variable with respect to these distributions allows for the general study of metric dimension, resolving sets, and efficient means for finding small resolving sets in some situations. Though there has not been as much work on metric dimension in this context as compared to deterministic graphs, there have been several significant contributions in this area concerning Erd\"os-R\'enyi random graphs~\cite{bollobas2012metric}, and random trees and forests~\cite{mitsche2015limiting}.

\subsection{Erd\"os-R\'enyi Random Graphs} For each $0\le p\le 1$, let $G_{n,p}$ denote a (simple) random graph with $n$ vertices obtained by including each of the possible ${n\choose 2}$ edges with probability $p$, independently of all other edges.

As $n \rightarrow \infty$, the set of $\lceil (3 \log n) / \log 2 \rceil$ highest degree vertices in $G_{n,1/2}$ suffices as a resolving set with high probability~\cite{babai1980random}. This upper bound on $\beta(G_{n,1/2})$ was originally used as part of a simple heuristic algorithm for canonically labeling graphs and determining whether or not two graphs are isomorphic~\cite{babai1980random}. More recently, focusing solely on adjacency information in $G_{n,p}$, this bound has been generalized for arbitrary values of $p$ to $\beta(G_{n,p}) \leq \frac{-3\ln(n)}{\ln(p^2+(1-p)^2)}$~\cite{TilLla19sbm,Till20}. The proof of this generalization does not rely on choosing a resolving set based on any particular property. In fact, any subset of nodes in $G_{n,p}$ of size at least $\frac{-3\ln(n)}{\ln(p^2+(1-p)^2)}$ is a resolving set with high probability for large $n$.

Detailed insight into the metric dimension of Erd\"os-R\'enyi graphs can be gained from the following result where $p$ is a function of $n$.

\begin{theorem}[Adapted from~\cite{bollobas2012metric}]
Let $d=(n-1)p$ be the expected degree.
Suppose that $$\log^5 n \ll d \leq n \left( 1-\frac{3\log(\log n)}{\log n} \right).$$
Let $i \geq 0$ be the largest integer such that $d^i = o(n)$ and let $c=c(n)=e^{d^{i+1}/n}$.
If $\beta_n$ denotes the metric dimension of $G_{n,p}$, then the following holds asymptotically almost surely.
\begin{displaymath}
 \beta_n =
 \begin{cases}
 \Theta(\log n) & \hbox{ if }c=\Theta(1),\\
 \Theta(c\log n) & \hbox{ if } c^{-1} = \Omega(d^i/n), \\
 \Theta(\frac{n \log n}{d^i}) & \hbox{ if } c^{-1} \ll d^i/n.
 \end{cases}
 \end{displaymath}
\end{theorem}

The regimes of $p$ described in this theorem produce a zig-zag pattern in $\beta(G_{n,p})$ as a function of $p$. Indeed, it can be shown that $\log_n(\beta(G_{n, n^{x-1}}))$, for $0<x<1$, approaches the function $f(x) = 1-x \lfloor 1/x \rfloor$ as $n\to\infty$ with high probability~\cite{bollobas2012metric}.

In a dense graph (i.e. with many edges and low path length entropy), consider picking a single vertex $v$ to add to a growing resolving set. This vertex defines an equivalence relation on the graph: two vertices are equivalent if they are the same distance away from $v$. In terms of distances and cardinalities, the equivalence classes are nearly the same---regardless of the vertex chosen. The ratio between the sizes of the two largest equivalence classes has great influence on the overall metric dimension. When this ratio is close to 1, picking a new vertex to add to the growing resolving set from the largest class will, on average, lead to more new equivalence classes than when the largest equivalence class contains many more vertices than the second largest. So, the overall metric dimension might be smaller the closer the two largest equivalence classes are in terms of size~\cite{bollobas2012metric}.

The zig--zag pattern observed in $\beta(G_{n,p})$ comes from how this ratio evolves in $G_{n, p}$ with decreasing $p$. Let $D_v(i)$ denote the set of vertices at a distance $i$ from a chosen vertex $v$. When $p=1$, the graph is complete and $|D_v(0)| = 1$ and $|D_v(1)| = n-1$ for all vertices $v$ so that $\beta(G_{n, 1}) = n-1$. As $p$ decreases, $|D_v(i)|$ for $i > 1$ increases. At first, this increase is faster for smaller values of $i$. Eventually, $|D_v(1)| \approx |D_v(2)|$ and the metric dimension is small. Decreasing $p$ further, the sizes of the two largest distance sets move away from one another and the ratio between their sizes increases. This pattern then repeats itself as the identities of the largest sets change~\cite{bollobas2012metric}.

Given two random variables, $X$ and $Y$, define the following measure of similarity between their distributions:
\begin{equation*}
    d(X,Y) := \sup_{h} \frac{E(h(X)) - E(h(Y))}{\sup_x |h(x)| + \sup_x |h'(x)|}
\end{equation*}
where the supremum is taken over all bounded test functions $h:\mathbb{R}\to\mathbb{R}$ with bounded derivative, and $E(\cdot)$ is used to denote expectation. The behavior of $\beta(G_{n,p})$ when $p$ is comparatively small is described by the following result. (The previous result addressed the case when $p$ is not as small.)
\begin{theorem}[\cite{mitsche2015limiting}]
Let $\beta_n$ denote the metric dimension of $G_{n,p}$.
\begin{itemize}
    \item[(i)] For $p=o(n^{-1})$, $\beta_n= n(1+o(1))$ asymptotically almost surely.
    \item[(ii)] For $p=c/n$ with $0 < c < 1$, the sequence of random variables
$$X_n = \frac{\beta_n - E(\beta_n)}{\sqrt{Var(\beta_n)}}$$
converges in distribution to a standard normal random variable $Z$ as $n \rightarrow \infty$, at a rate $d(X_n,Z) = O(n^{-1/2})$. 
Moreover, $E(\beta_n) = Cn(1+o(1))$ and $Var(\beta_n) = \Theta(n)$, where $C$ is an explicit constant that depends on $c$ only.
\end{itemize}
\end{theorem}

\subsection{Stochastic Block Model} The Stochastic Block Model (SBM) is a generative graph model used frequently to study networks with simple community structure. In its most basic form, the SBM has two main parameters, $C$ and $P$. $C$ is a partition of $n$ vertices into $c \geq 1$ disjoint communities $C_1, \dots, C_c$, and $P$ is a $(c \times c)$ symmetric matrix of adjacency probabilities. The communities are also sometimes defined stochastically using a probability vector of dimension $c$.

We say $G \sim SBM(n;C,P)$ when, for $u \in C_i$ and $v \in C_j$ with $u \neq v$, $\{u,v\}\in E$ with probability $P_{i,j}$, the entry in row-$i$ and column-$j$ of $P$, independently of all other pairs of nodes. Considered separately, each individual community is equivalent to an Erd\"os-R\'enyi random graph. This immediately suggests that approaches used to determine bounds on $\beta(G_{n,p})$ may prove valuable in this context as well. Inter-community adjacency probabilities, however, complicate the situation and must be dealt with carefully. For example, given $c=2$ and
\begin{displaymath}P = \begin{bmatrix}
o(1) & \frac{1}{2} \\
\frac{1}{2} & o(1)
\end{bmatrix}
\end{displaymath}
the communities of the resulting graph will be sparse and individual vertices may be difficult to distinguish without the help of vertices from both communities. Unfortunately, characterizing the precise interaction between vertices from different communities is not trivial, especially given the complicated dependencies between shortest path distances in such graphs.

It can be shown that the adjacency metric dimension of a graph, i.e. the metric dimension when only only neighbors of a node can be used to distinguish it from other nodes, serves as an upper bound on metric dimension~\cite{jannesari2012metric}. In particular, since the entries in the adjacency adjacency matrix of the SBM are by definition independent, a probabilistic upper bound on the metric dimension of these kinds of graphs can be established~\cite{TilLla19sbm,Till20}. Indeed, letting $\mathbb{P}(G;R)$ denote the probability that there are nodes in $G \sim SBM(n;C,P)$ with the same neighbors in $R$ where $R\subseteq\{1,\ldots,n\}$ contains $k_i$ nodes in community $i$, the first-moment method implies that
\begin{displaymath}
\mathbb{P}(G;R)\leq \sum_{1\leq i \leq j \leq c} |V_i||V_j|\prod_{\ell=1}^c P_{i,\ell}P_{j,\ell}+(1-P_{i,\ell})(1-P_{j,\ell})^{k_{\ell}}.
\end{displaymath}
This inequality serves as the basis of an effective, fast algorithm, for selecting vertices in $G$ such that $\sum_{1 \leq i \leq c} k_i$ is minimized and $\mathbb{P}(G;R)$ is less than a given threshold value. In essence, this algorithm provides an intelligent strategy for determining how the vertices of small resolving sets should be distributed across communities for any graph $G \sim SBM(n;C,P)$ with fixed parameters~\cite{TilLla19sbm,Till20}.

\subsection{Random Trees and Forests} 
\begin{figure}[h]
  \centering
  \includegraphics[scale=0.5]{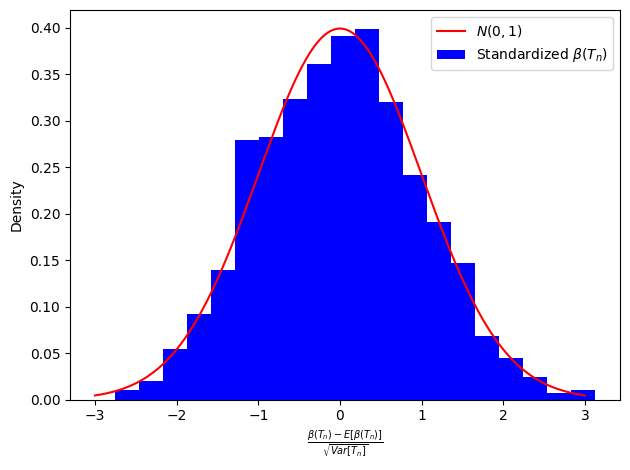}
  \vspace{-5mm}
  \caption[LoF entry]{Histogram of the metric dimensions of 1000 random trees on 2000 vertices (blue bars), along with the probability density function of a standard normal distribution (red curve). The Jarque-Bera test produces a $p$-value of $0.222$; in particular, the data carries little to no evidence against the null hypothesis that it comes from a standard Gaussian.}
  \label{fig:rand_tree_large_n}
\end{figure}

The metric dimension of a disconnected graph is, by definition, the sum of metric dimension of the graphs induced by its connected components. This is because the distance between any two vertices from different components is regarded as $\infty$. Accordingly, it is not surprising that the metric dimension of $F_n$, a forest on $n$ vertices chosen uniformly at random, has the same limiting distribution as that of $T_n$~\cite{mitsche2015limiting}, a tree on $n$ vertices also chosen uniformly at random. Furthermore, if $\beta_n=\beta(T_n)$ then, as $n \rightarrow \infty$, the random variables 
\begin{equation*}
    X_n = \frac{\beta_n - E(\beta_n)}{\sqrt{Var(\beta_n)}}
\end{equation*}
converge in distribution to a standard normal, where $E(\beta_n) = \mu n (1+o(1))$ and $Var(\beta_n) = \sigma^2 n(1+o(1))$, with $\mu \simeq 0.14076941$ and $\sigma^2 \simeq 0.063748151$ (see~\cref{fig:rand_tree_large_n}).

\section{Applications}

Small resolving sets are useful in a variety of situations. The direct analogy between metric dimension and trilateration in the plane makes potential applications regarding navigation~\cite{khuller1996landmarks} and location detection~\cite{slater1975leaves} in discrete space immediately apparent. Resolving sets have also been used as a means of comparing graphs. The classification of chemical compounds based on general chemical structure can be accomplished using resolving sets~\cite{chartrand2000resolvability,johnson1993structure,johnson1998browsable}. One method for quickly determining whether many, though not all, pairs of graphs are isomorphic relies on comparing vertex representations based on presumed resolving sets~\cite{babai1980random}. In this section, we examine several applications of resolving sets including as observers in detecting the source of a spread over a network~\cite{spinelli2016observer}, as a tool for detecting network motifs~\cite{hu2017detection}, and as the basis of a method for embedding DNA sequences in real space~\cite{TilLla18}.

\subsection{Source Localization}
\label{sec:source-loc}
A variety of transmission processes, such as information and disease, occur on networks. News and rumors circulate over social networks like Facebook and Twitter, and businesses take advantage of ``influencers'' to maximize the effect of marketing campaigns~\cite{booth2011mapping,katona2011network,kumar2012increasing,trusov2010determining}. Physical interaction networks are often used as a tool for studying the spread of diseases across communities~\cite{moore2000epidemics,newman2002spread,pastor2001epidemic,pinto2012locating}. In many contexts, it is valuable to locate the source of a transmission process, the node or set of nodes from which the spread began, for example to better understand the process and to decrease/increase transmission speed. Existing approaches to solving this problem include dynamic message passing~\cite{lokhov2014inferring}, time-reversal~\cite{shen2016locating}, and maximum likelihood estimators~\cite{pinto2012locating,shah2011rumors}.

Small resolving sets can also provide an elegant solution to source localization.
As a concrete example, suppose that a transmission process on the graph $G=(V,E)$ starts at an unknown source $s\in V$ at time $t_0$ and travels at unit speed across the edges.
Suppose we can pre-specify a set $R \subseteq V$ of observation nodes, such that we will observe the times $t_r$ at which the process first reaches each node $r\in R$.
A natural question is therefore, for which sets $R$ can we infer the source location $s$ from the observation times $\{t_r\}_{r\in R}$~\cite{pinto2012locating,zhang2016identification}?
The answer: exactly the sets $R$ which are resolving for $G$.
By definition of the process, we will have $t_r - t_0 = d(r,s)$, and thus the observation times uniquely identify $s$, for all possible sources $s$, exactly when $R$ is a resolving set.

While resolving sets present a promising solution to source localization, there are two non-trivial difficulties that arise in real-world transmission processes which require attention.
First, it is very unlikely that the start time $t_0$ of the process will be known.
Second, the speed $\ell_e$ at which information or diseases traverse a given edge $e\in E$, often called the edge length, is rarely deterministic or known precisely.
Instead, edge lengths are often modeled as randomly drawn from some known distribution. 

An unknown start time may be addressed by strengthening the notion of resolving sets to that of doubly resolving sets~\cite{caceres2007metric,spinelli2016observer}; see also \Cref{sec:doubly-resolving}.
We call $R \subseteq V$ \emph{doubly resolving} when for every pair of nodes $u,v \in V$ there is a pair $r, r'\in R$ such that $d(r, u) - d(r', u) \neq d(r, v) - d(r', v)$.
If $R$ is doubly resolving, it is also resolving.
Furthermore, one easily calculates $d(r, v) - d(r', v) = t_r - t_0 - t_{r'} + t_0 = t_r - t_{r'}$ for all nodes $v\in V$ and $r,r'\in R$.
Thus, as doubly resolving sets are based on relative and not absolute distances, the source can once again be located from the set $\{t_r\}_{r\in R}$ even if the start time $t_0$ is unknown.

Random edge lengths are more difficult to address and solutions depend on the details of the transmission process.
When the variance of $\ell_e$ is low relative to its mean, observation times will be close to their expected value, and resolving sets allow for exact solutions with high probability. When the variance of $\ell_e$ is high, however, observation times will carry very little information about the expected distances, especially for nodes $v$ at greater distance from the source $s$.
In this case, we can increase accuracy by adding nodes to $R$ using a path covering strategy such as truncated metric dimension (\Cref{sec:truncated}).
In this way, distances between observation nodes are small and thus noise accumulation is low, and we maintain resolvability of expected distances~\cite{spinelli2016observer}.

\subsection{Detecting Network Motifs}
\label{sec:motif}
A common tool in network science to compare graphs is via motifs, which are subgraphs appearing with higher than expected frequency.
Network motifs are believed to play important roles in the structure and underlying dynamics of networks in a variety of fields including social sciences~\cite{holland1974statistical,holland1975adaptation}, biology~\cite{eichenberger2004program,gosak2018network,lee2002transcriptional,shen2002network}, chemistry~\cite{hahnke2010pharmacophore,weininger1988smiles}, and data mining more generally~\cite{conte2004thirty,von1988pattern,washio2003state}. By finding and analyzing these motifs, researchers gain insight into the functional properties of different systems. The problem of discovering important subgraphs in a large network, however, poses significant computational challenges: the subgraph isomorphism problem, determining whether or not a given graph occurs as a subgraph in a larger graph, is NP-complete~\cite{cook1971complexity,wegener2005complexity}.
On the other hand, graph isomorphism, the special case of subgraph isomorphism when the graphs have the same size, is believed to be an easier problem~\cite{schoning1988graph,crasmaru2004protocol}, especially given the recent quasi-polynomial time algorithm~\cite{babai2016graph}. A natural algorithm to count $k$-node motifs is therefore to enumerate subgraphs of size $k$ and test whether they are isomorphic to the given motif.
Resolving sets and metric dimension have been used as the foundation of tools to solve graph isomorphism~\cite{babai1980random}, and in the manner above, for motif detection~\cite{hu2017detection}.

One technique to solve graph isomorphism is through a \emph{canonical labeling}, a way to assign unique labels to nodes which is invariant under graph isomorphism.
Since resolving sets provide a unique label for each node in a graph, they can serve as the basis of such labelings.
In particular, suppose that one could uniquely identify a resolving set $R$ given a graph $G=(V,E)$, and furthermore, uniquely identify an ordering $R = \{r_1,\ldots,r_{|R|}\}$.
Then assigning label $d(u|R)\in\mathbb R^{|R|}$ to each node $u\in V$ will give a canonical labeling, as both $R$ and the ordering is uniquely determined by $G$, and given $R$ and this ordering, the distance vectors are unique as $R$ is resolving.

Using this general approach, the following quadratic-time canonical labeling algorithm provably solves graph isomorphism with high probability for Erd\"{o}s-R\'{e}nyi random graphs $G_{n,p}$, that is, graphs on $n$ vertices such that each edge appears with independent probability $p$~\cite{babai1980random}.
Consider a graph $G=(V,E)$ and let $n=|V|$ and $R=\{r_1, \dots, r_{|R|}\}$ be the set of the $\lceil (3 \log(n)) / \log(2) \rceil$ highest degree vertices in $V$.\footnote{In fact, any set of vertices of this size will suffice as a resolving set with high probability, not just vertices of high degree. This fact is directly related to the fact that the degree distribution of $G_{n,p}$ is Binomial, and therefore concentrated around its mean~\cite{arratia1990erdos,cramer1938nouveau}.}
The algorithm labels each vertex in $G$ with the set of nodes in $R$ which are adjacent.
(This step of the algorithm leverages the fact that the diameter of $G_{n,p}$ is 2 with high probability as $n$ increases, so the set of adjacent $R$ nodes is equivalent to the distance vector $d(\cdot|R)$.)
This labeling is canonical, i.e., invariant under isomorphism, as long as no pair of vertices in $R$ have the same degree, and the labels are unique. Under $G_{n,1/2}$, Erd\"{o}s-R\'{e}nyi random graphs with adjacency probability $\frac{1}{2}$, the probability that this algorithm succeeds in finding a canonical labeling is at least $1-\sqrt[7]{1/n}$ for large $n$~\cite{babai1980random}. This technique indirectly shows that $\beta(G_{n,1/2}) \leq \lceil (3 \log(n)) / \log(2) \rceil$; more generally, it can be shown that $\beta(G_{n,p}) \leq \frac{-3\ln(n)}{\ln(p^2+(1-p)^2)}$~\cite{TilLla19sbm,Till20}.

Another recent algorithm uses the above technique more explicitly, by directly computing a canonical resolving set and ordering~\cite{hu2017detection}. 
The approach is as follows.
Enumerate all resolving sets $R$ of size $\beta(G)$ and permutations $\pi$, and compute the corresponding adjacency matrix $A_{R,\pi}$ for each.
In $A_{R,\pi}$, the vertices of $G$ are ordered lexicographically with respect to their distance vectors $d(\cdot|R)$ ordered by $\pi$.
Now take $R$ and $\pi$ such that $A_{R,\pi}$ is lexicographically first (after flattening) among all such choices.
Given this choice of $R$ and $\pi$, labeling the vertices of $G$ by their distance vector representations gives a canonical labeling.

This approach can be time consuming, as it requires computing (an upper bound for) $\beta(G)$ and enumerating all resolving sets $R$ of this size and all orderings on these sets.
One can improve performance slightly by ignoring the relative ordering of twin vertices (\S~\ref{sec:twin_nodes})~\cite{hu2017detection}. Tests on a wide variety of graphs show that this method is effective but somewhat slower than the graph isomorphism tool provided with nauty~\cite{mckay2014practical}. However, this method is faster than nauty on multi-dimensional mesh graphs, and was the core algorithm in a tool for identifying and counting statistically significant subgraphs in the transcriptional regulation networks of \emph{Saccharomyces cerevisiae} (yeast) and \emph{Escherichia coli} (\emph{E. coli})~\cite{hu2017detection}.

\subsection{Embedding Biological Sequence Data}
High-throughput sequencing technologies have enabled biologists to collect a wealth of DNA, RNA, and amino acid sequence data. The abundance of this information makes computational analysis methods, including those based on machine learning algorithms, indispensable.
The majority of these methods, however, cannot directly learn from symbolic data, like biological sequences, and deal instead with numeric vectors.
Methods to embed symbolic sequences into real vector spaces are thus an important pre-processing step~\cite{TilLla18}.
Low-dimensional embeddings are especially useful, both to reduce the computational cost of learning algorithms and to avoid overfitting.

Consider the task of embedding a sequence of length $\ell$, composed of symbols from an alphabet of size $a$, into a real vector space.
For example, DNA and RNA sequences have $a=4$, while amino acids are composed of $a=20$ possible symbols.
A na\"ive approach to this embedding is the so-called ``one-hot encoding'', which simply generates an indicator vector for each of the $a^\ell$ possible sequences of length $\ell$.
This approach is untenable for most biological sequence data, where $\ell$ can be quite large.
One lower-dimensional approach to embed such sequences uses binary representations, which indicate the presence or absence of each character in the alphabet at each position in the sequence, thus requiring $a\cdot \ell$ dimensions~\cite{cai2003support}.
Another common approach uses $k$-mer count vectors, which count the number of times that every possible contiguous subsequence of length $k$ occurs in the larger sequence using a sliding window~\cite{leslie2002spectrum}.
Unlike one-hot encodings and binary representations, $k$-mer count vectors are not guaranteed to produce an injective embedding.
These $k$-mer count vectors are $a^k$-dimensional; typically one chooses $k \ll \ell$.

Resolving sets can be applied to this problem a well.
Often the domain of interest has some natural distance metric between sequences.
For biological sequences and several other domains, a natural choice is Hamming distance, which simply counts the number of indices in which the two sequences differ.
The Hamming distance induces the Hamming graph $H_{\ell,a}$ on length-$\ell$ sequences from $a$ symbols, where there are edges between sequences which differ in only one entry, and hence the path length between two sequences is their Hamming distance (\Cref{sec:hamming}).
Given a resolving set $R$ on $H_{\ell,a}$, each sequence can be uniquely represented by its vector of distances to the elements of $R$.
These distance vectors are therefore an injective $|R|$-dimensional embedding, and the metric dimension of $H_{\ell,a}$ gives the smallest possible embedding dimension for this approach.

Embeddings based on resolving sets of $H_{3,4}$ were used as features to classify DNA sequences of length 20 as being centered at intron-exon boundaries, or not, in the fruit fly genome~\cite{TilLla18}.
In this study, the resolving set embedding outperforms $k$-mer count vector and binary representation based features with respect to accuracy and specificity, and is competitive with features based on other state-of-the-art embedding techniques like Node2Vec~\cite{grover2016node2vec} and multidimensional scaling~\cite{Krz00}.\footnote{For this particular task, $k$-mer count vectors are not well suited. As positive examples consist of half intronic and half exonic DNA, we might expect the location of $k$-mers within the larger sequence to matter a great deal, yet $k$-mer count vectors do not directly encode this information.}
Resolving set embeddings are also generally more compact than those based on $k$-mer count vectors or binary representations.
For example, part of the genome of the Dengue virus codes for a protease that targets octapeptides, amino acid sequences of length 8, in human cells.
While the space of all octapeptides is large, consisting of $25.6$ billion sequences, no more than $82$ are required for a resolving set.
Based on such a set, $H_{8,20}$ may be embedded in $\mathbb{R}^{82}$. 
In comparison, $3$-mer count vectors use $8,000$ dimensions while a binary vector representation requires 160 dimensions~\cite{TilLla18}.

\section{Related Concepts}

There are a wide variety of concepts closely linked to metric dimension.
Some strengthen or alter the constraints placed on the identifiability of nodes while others are extensions of the concept itself.
In this section we will briefly define several such notions.
A survey by Chartrand and Zhang~\cite{chartrand2003theory} contains more complete characterizations of several of these concepts as well as information on concepts not mentioned here.
There is also a body of work on conditional resolvability which focuses on resolving sets that have some additional property.
For example, one can consider the smallest resolving set of a graph that induces a connected subgraph, or that is also an independent set.
For these conditional variants, we direct the reader to a survey by Saenpholphat and Zhang~\cite{saenpholphat2004conditional}.

\subsection{Doubly Resolving Sets}
\label{sec:doubly-resolving}

As discussed in \Cref{sec:source-loc}, while one can uniquely identify all nodes in a graph based on distances to a resolving set, this identification can fail if one only knows distances up to an additive constant.
Doubly resolving sets address this shortcoming and have proven useful in identifying the source of a spread in a network~\cite{spinelli2016observer} and in determining bounds on the metric dimension of Cartesian products of graphs~\cite{caceres2007metric}. In particular, a set $R \subseteq V$ is called \emph{doubly resolving} if, for every pair of nodes $u,v \in V$, there is a pair $r, r' \in R$ such that $d(r,u) - d(r',u) \neq d(r,v) - d(r',v)$. Such sets differentiate nodes based on relative as opposed to absolute distances.

For example, consider the path graph $P_n$ with nodes labeled consecutively from $1$ to $n$ and resolving set $\{1\}$. Suppose a signal is sent along the path from some node $i$ at arbitrary time $t$. Traversing each edge in one time unit, this signal reaches node $1$ at time $(t+i-1)$. Since $t$ is unknown, node $i$ cannot be distinguished as the source of the signal from this information alone. 
The set $\{1, n\}$, on the other hand, is doubly resolving on $P_n$. Now a signal sent from $i$ at time $t$ will arrive at node $1$ at time $(t+i-1)$ and at node $n$ at time $(t+n-i)$. So, for any node $j \neq i$, $d(i, 1) - d(i, n) = (t+i-1)-(t+n-i) = 2i-n-1 \neq 2j-n-1 = (t+j-1)-(t+n-j) = d(j, 1) - d(j, n)$ and the source of the signal can be uniquely determined. 

\subsection{Strong Metric Dimension}
\label{sec:strong-metric-dim}
Since all vertices in a graph $G=(V,E)$ are distinguished based on distances to a resolving set $R \subseteq V$, it is tempting to think that $G$ may be reconstructed using $R$.
These distances can only recover shortest paths, however, and thus any edge in $E$ that is not part of a unique shortest path from any $v \in V$ to any $r \in R$ could be excluded in such a reconstruction.
For example, consider the cycle $C_6$ with the minimum resolving set $R=\{1,3\}$. Notice $R$ remains a minimum resolving set if the edge $(2,5)$ is added. Furthermore, for each $v \in V$, $d(v|R)$ is the same with or without this extra edge. As a result, $R$ is not enough to guarantee a faithful reconstruction of $C_6$.
Notably, any resolving set is enough to reconstruct a tree, as the only path between any pair of vertices is the shortest path.

A set $S \subseteq V$ is said to \emph{strongly resolve} $G$ if for every $u,v \in V$ there is a vertex $s \in S$ such that $u$ lies on a shortest path from $s$ to $v$ or $v$ lies on a shortest path from $s$ to $u$~\cite{sebHo2004metric}. The size of smallest possible strongly resolving sets on a graph is its strong metric dimension. By definition, every edge in $E$ must be accounted for by a shortest path distance $d(v, s)$ for some $v \in V$ and $s \in S$. This allows $G$ to be reconstructed exactly based on a strong resolving set~\cite{sebHo2004metric}. For a survey of results and approximation methods related to strong metric dimension see~\cite{kratica2016strong}.

\subsection{Multilateration}
\label{sec:multilateration}
The definition of metric dimension depends heavily on graph structure and a notion of edge distances. Yet there are pairwise distance matrices that do not correspond to a graph or metric space. Consider the following matrix:
\begin{displaymath}M = \begin{blockarray}{cccc}
& A & B & C \\
\begin{block}{c(ccc)}
  A & 0 & 10 & 100 \\
  B & \infty & 0 & 10 \\
  C & \infty & \infty & 0 \\
\end{block}
\end{blockarray}
\end{displaymath}

\noindent Notice that the points in $M$ do not abide by the triangle inequality. While we can circumvent the fact that values in this matrix are not symmetric by using directed edges, violating the triangle inequality would require redefining the distance between pairs of nodes.
Fortunately, metric dimension makes no use of the actual distance between two vertices beyond checking equality, and thus one can imagine a relaxation using a more general distance function.
In fact, one could even allow the entries of $M$ to come from an arbitrary set other than the reals.

This more general problem on arbitrary matrices is called \emph{multilateration}~\cite{TilLla18}.
Let $I$ be a set of items associated with the rows of a matrix $M$ and $F$ be a set of functions over $I$ associated with columns such that $M(i,f) = f(i)$.
Analogous to metric dimension, the goal of multilateration is to determine a resolving set $R \subseteq F$ of minimum size such that the vectors $(r(i))_{r\in R}$ are unique for all $i \in I$. One can equivalently think of $R$ as a set of columns of $M$ such that the row vectors of the induced submatrix are unique. Borrowing notation from metric dimension, we set $\beta(M)=|R|$ for a set $R$ of minimum size.
The entries $M(i,f)$ need not be numeric; one only needs a notion of equivalence on elements of $f(I)$ for each $f\in F$, i.e., for values in the same column of $M$.
In the context of graphs, multilateration is equivalent to metric dimension: if $G$ is a graph with pairwise distance matrix $D$, then we have $\beta(G) = \beta(D)$.

\subsection{Truncated Metric Dimension}
\label{sec:truncated}
In some scenarios, complete distance information of a network is unavailable.
In particular, perhaps only distances below a certain threshold are available, perhaps because collecting long-distance information is costly or prone to an excessive amount of noise.
In such cases it may not be possible to determine the metric dimension of the full graph. Instead, given a graph $G=(V,E)$ and a maximum distinguishable distance $k$, let $d_k(u,v)=\min\{d(u,v),k+1\}$ be the $k$-truncated distance between $u,v \in V$ and let $D_k$ be the $k$-truncated distance matrix of $G$.
Then the $k$-truncated metric dimension of $G$ is defined as $\beta_k(G) = \beta(D_k)$, where $\beta(D_k)$ is defined as in multilateration (\Cref{sec:multilateration})~\cite{TilLla19dmini,Till20}. This notion is a generalization of the concept of adjacency metric dimension~\cite{jannesari2012metric} where $k=1$ and vertices are distinguished by their neighbors in a resolving set.

Beyond settings with restricted distance information, truncated metric dimension can also be an effective tool for studying metric dimension.
For instance, it can be shown that for all graphs $G$ and all $k\geq 1$, we have $\beta(D_k) \geq \beta(D_{k+1})$~\cite{TilLla19dmini,Till20}; as $\beta(G) = \beta(D_{n-1})$, this means truncated metric dimension can give upper bounds on $\beta(G)$. As one application, an asymptotically tight upper bound on the metric dimension of the Erd\"os-R\'enyi random graph $G_{n,p}$ can be determined by focusing on $\beta(D_1)$, i.e., on adjacency information alone~\cite{TilLla19sbm,Till20}. While $G_{n,p}$ contains dependencies between shortest path lengths, adjacencies are independent by definition, making the $1$-truncated metric dimension of these graphs far easier to characterize than their standard metric dimension (\Cref{sec:motif}).

\subsection{Resolving Number, Upper Dimension, Random $k$-dimensionality}
Given a resolving set $R$ of $G=(V,E)$, it is natural to consider removing elements of $R$ while keeping the set resolving.
Let $\mathcal R(G)$ be the set of all resolving sets which cannot be made smaller in this way, i.e., for which no proper subset $S \subset R$ is also resolving.
Resolving sets $R \in \mathcal R(G)$, while minimal in the sense of set inclusion, are not guaranteed to be minimal in the sense that $|R|=\beta(G)$, the minimum possible size of any resolving set of $G$.
For example, on the path $P_6$ on vertices $\{1,\ldots,6\}$, the set $R=\{3,4\}$ is in $\mathcal R(P_6)$, but $\beta(P_6)=1$, achieved by $\{1\}$. As resolving sets of cardinality $\beta(G)$ are also elements of $\mathcal R(G)$, however, we do have $\beta(G) = \min_{R \in \mathcal R(G)} |R|$.
The size of largest set-inclusion-minimal resolving set, on the other hand, is called the \emph{upper dimension} $\mathrm{dim}^+(G)$, given by $\mathrm{dim}^+(G) = \max_{R \in \mathcal R(G)} |R|$.

To generate sets in $\mathcal R(G)$, one could start with $R = V$, and iteratively remove vertices while keeping $R$ resolving.
In some cases, starting with a smaller set is guaranteed to succeed as well.
The \emph{resolving number} of a graph, denoted $\mathrm{res}(G)$, is the smallest integer such that all subsets $S \subset V$ with $|S| = \mathrm{res}(G)$ are resolving sets.
Combined with the above inequalities, we therefore have $\beta(G) \leq dim^+(G) \leq \mathrm{res}(G) \leq n-1$~\cite{chartrand2000resolvabilityUpperDim}.
When $k = \beta(G) = \mathrm{res}(G)$, every $R \in \mathcal R(G)$ with $|R| = k$ is resolving and is of minimum size.
In this case $G$ is said to be \emph{randomly $k$-dimensional}.
The only known randomly $k$-dimensional graphs are $K_{k+1}$ and odd cycles $C_n$ with $n \geq 3$ and $k=2$. It is an open question as to whether or not other types of randomly $k$-dimensional graphs exist~\cite{chartrand2000chromatic}.

\section{Conclusion}

Intuitively, metric dimension is a very simple idea. Its close relation to GPS and trilateration in continuous space make applications concerning locating nodes of graphs immediately clear.
While determining the exact metric dimension of general graphs is an NP-complete problem, a number of approximation methods exist including the ICH algorithm.
ICH gives small resolving sets guaranteed to be close to optimal but has cubic run time, which can be impractical for large networks.
Fortunately, a great deal is known with respect to exact formulae, asymptotic behavior, and bounds for an array of graph families.
Recently, the metric dimension of certain random graph models like random trees and forests, Erd\"os-R\'enyi random graphs, and graphs generated via the SBM have been characterized.
Further work in this direction may lead to a more complete understanding of metric dimension and of its behavior on real-world networks.

This understanding is a critical step toward practical application in different settings.
These applications include identifying the source of a spread in a network, detecting network motifs, and embedding symbolic data in real space.
Some of these applications make use of concepts related to metric dimension allowing slightly different constraints, such as doubly resolving sets and multilateration.
Possible directions for future work abound, both in deepening our theoretical understanding of metric dimension and related concepts, and in applying these concepts in practice. 

\section*{Acknowledgements}
This research was partially funded by NSF ISS grant No. 1836914. The authors acknowledge the BioFrontiers Computing Core at the University of Colorado--Boulder for providing High-Performance Computing resources (fund\-ed by the NIH grant No. 1S10OD012300), supported by BioFrontiers IT group.

\bibliographystyle{siam}
\bibliography{biblio.bib}
\end{document}